\documentclass{gtart}

\def\ifplaintex{\expandafter\ifx\csname documentclass\endcsname\relax}

\def\gtp{{\mathsurround=0pt\it $\cal G\mskip-2mu$eometry \&\ 
$\cal T\!\!$opology $\cal P\!$ublications}}  

\def\recd{{\small Received:\qua\receiveddate\ifx\reviseddate\relax
\else\qquad Revised:\qua\reviseddate\fi\par}} 

\def\lognumber#1{\def\thelognumber{#1}}
\def\volumenumber#1{\def\thevolumenumber{#1}}
\def\volumeyear#1{\def\thevolumeyear{#1}}
\def\papernumber#1{\def\thepapernumber{#1}}
\def\pagenumbers#1#2{\def\startpage{#1}\def\finishpage{#2}}
\def\published#1{\def\publishdate{#1}}

\def\received#1{\def\receiveddate{#1}}
\def\revised#1{\def\reviseddate{#1}}
\def\accepted#1{\def\accepteddate{#1}}

\long\def\asciiabstract#1{\long\def\theasciiabstract{#1}}

\let\\\par\let\thelognumber\relax\let\thevolumenumber\relax
\let\thepapernumber\relax\let\thevolumeyear\relax\let\startpage\relax
\let\finishpage\relax\let\publishdate\relax\let\receiveddate\relax
\let\reviseddate\relax\let\accepteddate\relax\let\theasciititle\relax
\let\theasciiauthors\relax
\let\theasciiabstract\relax

\let\theasciiemail\relax

\ifplaintex
\font\logobig=cmssbx10 scaled 3836
\font\logomed=cmssbx10 scaled 2557
\else
\font\logobig=cmssbx10 scaled 4200
\font\logomed=cmssbx10 scaled 2800
\fi

\long\def\makeagttitle{   
\count0=\startpage
\agt\hfill      
\hbox to 45truept{\vbox to 0pt{\vglue -13truept{\logomed A\kern -.37em{\logobig 
T}\kern -.38em G}\vss}\hss}
\break
{\small Volume \thevolumenumber\ (\thevolumeyear)
\startpage--\finishpage\nl
Published: \publishdate}

\vglue .25truein

{\parskip=0pt\leftskip 0pt plus
1fil\def\\{\par\smallskip}{\Large\bf\thetitle}\par\medskip} \vglue
0.05truein

{\parskip=0pt\leftskip 0pt plus 1fil\def\\{\par}{\sc\theauthors}
\par\medskip}%
 
\vglue 0.03truein 

{\small\leftskip 25truept\rightskip 25truept{\bf Abstract}\stdspace\theabstract

{\bf AMS Classification}\stdspace\theprimaryclass
\ifx\thesecondaryclass\relax\else; \thesecondaryclass\fi\par
{\bf Keywords}\stdspace \thekeywords\par}\vglue 7truept

}   

\ifplaintex
\font\phead=cmsl9 scaled 950
\font\pnum=cmbx10 scaled 913
\font\pfoot=cmsl9 scaled 950
\headline{\vbox to 0pt{\vskip -4.5mm\line{\small\phead\ifnum
\count0=\startpage ISSN 1472-2739 (on-line) 1472-2747 (printed)
\hfill {\pnum\folio}\else\ifodd\count0\def\\{ }%
\ifx\theshorttitle\relax\thetitle\else\theshorttitle\fi\hfill{\pnum\folio}
\else\def\\{ and }{\pnum\folio}\hfill\ifx\theshortauthors\relax\theauthors
\else\theshortauthors\fi\fi\fi}\vss}}
\footline{\vbox to 0pt{\vglue 0mm\line{\small\pfoot\ifnum\count0=\startpage
\copyright\ \gtp\hfill\else
\agt, Volume \thevolumenumber\ (\thevolumeyear)\hfill\fi}\vss}}
\else
\font=cmsl9 scaled 1050
\font\lnum=cmbx10 
\font=cmsl9 scaled 1050
\makeatletter
\def\@oddhead{{\small\lhead\ifnum\count0=\startpage ISSN 1472-2739 
(on-line) 1472-2747 (printed)\hfill {\lnum\number\count0}\else\ifodd\count0
\def\\{ }\ifx\theshorttitle\relax \thetitle \else\theshorttitle\fi\hfill
{\lnum\number\count0}\else\def\\{ and }{\lnum\number\count0}
\hfill\ifx\theshortauthors\relax 
\theauthors\else\theshortauthors\fi\fi\fi}}\def\@evenhead{\@oddhead}
\def\@oddfoot{\small\lfoot\ifnum\count0=\startpage\copyright\ \gtp\hfill\else
\agt, Volume \thevolumenumber\ (\thevolumeyear)\hfill\fi}
\def\@evenfoot{\@oddfoot}
\makeatother
\fi
\let\maketitlepage\makeagttitle
\let\makeshorttitle\maketitlepage
\let\maketitle\maketitlepage

\newwrite\gtoutfile
\long\gdef\makeheadfile{  
{\def\\{, }\def\s{ }
\immediate\openout\gtoutfile head.xxx
\immediate\write\gtoutfile{To: math@arxiv.org}
\immediate\write\gtoutfile{Subject: put OR rep NNNNN:ppppp}
\immediate\write\gtoutfile{--text follows this line--}
\immediate\write\gtoutfile{Proxy-for: \ifx\theasciiauthors\relax
\theauthors\else\theasciiauthors\fi\s<\ifx\theasciiemail\relax\theemail\else\theasciiemail\fi>}
\immediate\write\gtoutfile{\noexpand\\}
\immediate\write\gtoutfile{Authors: \ifx\theasciiauthors\relax
\theauthors\else\theasciiauthors\fi}
{\def\\{ }\immediate\write\gtoutfile{Title: \ifx\theasciititle\relax
\thetitle\else\theasciititle\fi}}
\immediate\write\gtoutfile{Subj-class: GT or SG, GR etc}
\immediate\write\gtoutfile{MSC-class: \theprimaryclass\ifx\thesecondaryclass\relax\else, \thesecondaryclass\fi}
\immediate\write\gtoutfile{Journal-ref: Algebr. Geom. Topol. \thevolumenumber\s
(\thevolumeyear) \startpage-\finishpage}
\immediate\write\gtoutfile{Comments: Published by Algebraic and
Geometric Topology at}
\immediate\write\gtoutfile{\s\s\s  http://www.maths.warwick.ac.uk/agt/AGTVol\thevolumenumber/agt-\thevolumenumber-\thepapernumber.abs.html}
\immediate\write\gtoutfile{\noexpand\\}
\immediate\write\gtoutfile{}
\ifx\theasciiabstract\relax
\immediate\write\gtoutfile{\theabstract}\else
\immediate\write\gtoutfile{\theasciiabstract}\fi
\immediate\write\gtoutfile{}
\immediate\write\gtoutfile{\noexpand\\}
\immediate\write\gtoutfile{}
\immediate\closeout\gtoutfile}}  

\def\maketitlepage{\makeagttitle\makeheadfile}
\let\makeshorttitle\maketitlepage
\let\maketitle\maketitlepage

\def\ifplaintex{\expandafter\ifx\csname documentclass\endcsname\relax}

\def\gtp{{\mathsurround=0pt\it $\cal G\mskip-2mu$eometry \&\ 
$\cal T\!\!$opology $\cal P\!$ublications}}  

\def\recd{{\small Received:\qua\receiveddate\ifx\reviseddate\relax
\else\qquad Revised:\qua\reviseddate\fi\par}} 

\def\lognumber#1{\def\thelognumber{#1}}
\def\volumenumber#1{\def\thevolumenumber{#1}}
\def\volumeyear#1{\def\thevolumeyear{#1}}
\def\papernumber#1{\def\thepapernumber{#1}}
\def\pagenumbers#1#2{\def\startpage{#1}\def\finishpage{#2}}
\def\published#1{\def\publishdate{#1}}

\def\received#1{\def\receiveddate{#1}}
\def\revised#1{\def\reviseddate{#1}}
\def\accepted#1{\def\accepteddate{#1}}

\long\def\asciiabstract#1{\long\def\theasciiabstract{#1}}

\let\\\par\let\thelognumber\relax\let\thevolumenumber\relax
\let\thepapernumber\relax\let\thevolumeyear\relax\let\startpage\relax
\let\finishpage\relax\let\publishdate\relax\let\receiveddate\relax
\let\reviseddate\relax\let\accepteddate\relax\let\theasciititle\relax
\let\theasciiauthors\relax
\let\theasciiabstract\relax

\let\theasciiemail\relax

\ifplaintex
\font\logobig=cmssbx10 scaled 3836
\font\logomed=cmssbx10 scaled 2557
\else
\font\logobig=cmssbx10 scaled 4200
\font\logomed=cmssbx10 scaled 2800
\fi

\long\def\makeagttitle{   
\count0=\startpage
\agt\hfill      
\hbox to 45truept{\vbox to 0pt{\vglue -13truept{\logomed A\kern -.37em{\logobig 
T}\kern -.38em G}\vss}\hss}
\break
{\small Volume \thevolumenumber\ (\thevolumeyear)
\startpage--\finishpage\nl
Published: \publishdate}

\vglue .25truein

{\parskip=0pt\leftskip 0pt plus
1fil\def\\{\par\smallskip}{\Large\bf\thetitle}\par\medskip} \vglue
0.05truein

{\parskip=0pt\leftskip 0pt plus 1fil\def\\{\par}{\sc\theauthors}
\par\medskip}%
 
\vglue 0.03truein 

{\small\leftskip 25truept\rightskip 25truept{\bf Abstract}\stdspace\theabstract

{\bf AMS Classification}\stdspace\theprimaryclass
\ifx\thesecondaryclass\relax\else; \thesecondaryclass\fi\par
{\bf Keywords}\stdspace \thekeywords\par}\vglue 7truept

}   

\ifplaintex
\font\phead=cmsl9 scaled 950
\font\pnum=cmbx10 scaled 913
\font\pfoot=cmsl9 scaled 950
\headline{\vbox to 0pt{\vskip -4.5mm\line{\small\phead\ifnum
\count0=\startpage ISSN 1472-2739 (on-line) 1472-2747 (printed)
\hfill {\pnum\folio}\else\ifodd\count0\def\\{ }%
\ifx\theshorttitle\relax\thetitle\else\theshorttitle\fi\hfill{\pnum\folio}
\else\def\\{ and }{\pnum\folio}\hfill\ifx\theshortauthors\relax\theauthors
\else\theshortauthors\fi\fi\fi}\vss}}
\footline{\vbox to 0pt{\vglue 0mm\line{\small\pfoot\ifnum\count0=\startpage
\copyright\ \gtp\hfill\else
\agt, Volume \thevolumenumber\ (\thevolumeyear)\hfill\fi}\vss}}
\else
\font=cmsl9 scaled 1050
\font\lnum=cmbx10 
\font=cmsl9 scaled 1050
\makeatletter
\def\@oddhead{{\small\lhead\ifnum\count0=\startpage ISSN 1472-2739 
(on-line) 1472-2747 (printed)\hfill {\lnum\number\count0}\else\ifodd\count0
\def\\{ }\ifx\theshorttitle\relax \thetitle \else\theshorttitle\fi\hfill
{\lnum\number\count0}\else\def\\{ and }{\lnum\number\count0}
\hfill\ifx\theshortauthors\relax 
\theauthors\else\theshortauthors\fi\fi\fi}}\def\@evenhead{\@oddhead}
\def\@oddfoot{\small\lfoot\ifnum\count0=\startpage\copyright\ \gtp\hfill\else
\agt, Volume \thevolumenumber\ (\thevolumeyear)\hfill\fi}
\def\@evenfoot{\@oddfoot}
\makeatother
\fi
\let\maketitlepage\makeagttitle
\let\makeshorttitle\maketitlepage
\let\maketitle\maketitlepage

\newwrite\gtoutfile
\long\gdef\makeheadfile{  
{\def\\{, }\def\s{ }
\immediate\openout\gtoutfile head.xxx
\immediate\write\gtoutfile{To: math@arxiv.org}
\immediate\write\gtoutfile{Subject: put OR rep NNNNN:ppppp}
\immediate\write\gtoutfile{--text follows this line--}
\immediate\write\gtoutfile{Proxy-for: \ifx\theasciiauthors\relax
\theauthors\else\theasciiauthors\fi\s<\ifx\theasciiemail\relax\theemail\else\theasciiemail\fi>}
\immediate\write\gtoutfile{\noexpand\\}
\immediate\write\gtoutfile{Authors: \ifx\theasciiauthors\relax
\theauthors\else\theasciiauthors\fi}
{\def\\{ }\immediate\write\gtoutfile{Title: \ifx\theasciititle\relax
\thetitle\else\theasciititle\fi}}
\immediate\write\gtoutfile{Subj-class: GT or SG, GR etc}
\immediate\write\gtoutfile{MSC-class: \theprimaryclass\ifx\thesecondaryclass\relax\else, \thesecondaryclass\fi}
\immediate\write\gtoutfile{Journal-ref: Algebr. Geom. Topol. \thevolumenumber\s
(\thevolumeyear) \startpage-\finishpage}
\immediate\write\gtoutfile{Comments: Published by Algebraic and
Geometric Topology at}
\immediate\write\gtoutfile{\s\s\s  http://www.maths.warwick.ac.uk/agt/AGTVol\thevolumenumber/agt-\thevolumenumber-\thepapernumber.abs.html}
\immediate\write\gtoutfile{\noexpand\\}
\immediate\write\gtoutfile{}
\ifx\theasciiabstract\relax
\immediate\write\gtoutfile{\theabstract}\else
\immediate\write\gtoutfile{\theasciiabstract}\fi
\immediate\write\gtoutfile{}
\immediate\write\gtoutfile{\noexpand\\}
\immediate\write\gtoutfile{}
\immediate\closeout\gtoutfile}}  

\def\maketitlepage{\makeagttitle\makeheadfile}
\let\makeshorttitle\maketitlepage
\let\maketitle\maketitlepage

\def\ifplaintex{\expandafter\ifx\csname documentclass\endcsname\relax}

\def\gtp{{\mathsurround=0pt\it $\cal G\mskip-2mu$eometry \&\ 
$\cal T\!\!$opology $\cal P\!$ublications}}  

\def\recd{{\small Received:\qua\receiveddate\ifx\reviseddate\relax
\else\qquad Revised:\qua\reviseddate\fi\par}} 

\def\lognumber#1{\def\thelognumber{#1}}
\def\volumenumber#1{\def\thevolumenumber{#1}}
\def\volumeyear#1{\def\thevolumeyear{#1}}
\def\papernumber#1{\def\thepapernumber{#1}}
\def\pagenumbers#1#2{\def\startpage{#1}\def\finishpage{#2}}
\def\published#1{\def\publishdate{#1}}

\def\received#1{\def\receiveddate{#1}}
\def\revised#1{\def\reviseddate{#1}}
\def\accepted#1{\def\accepteddate{#1}}

\long\def\asciiabstract#1{\long\def\theasciiabstract{#1}}

\let\\\par\let\thelognumber\relax\let\thevolumenumber\relax
\let\thepapernumber\relax\let\thevolumeyear\relax\let\startpage\relax
\let\finishpage\relax\let\publishdate\relax\let\receiveddate\relax
\let\reviseddate\relax\let\accepteddate\relax\let\theasciititle\relax
\let\theasciiauthors\relax
\let\theasciiabstract\relax

\let\theasciiemail\relax

\ifplaintex
\font\logobig=cmssbx10 scaled 3836
\font\logomed=cmssbx10 scaled 2557
\else
\font\logobig=cmssbx10 scaled 4200
\font\logomed=cmssbx10 scaled 2800
\fi

\long\def\makeagttitle{   
\count0=\startpage
\agt\hfill      
\hbox to 45truept{\vbox to 0pt{\vglue -13truept{\logomed A\kern -.37em{\logobig 
T}\kern -.38em G}\vss}\hss}
\break
{\small Volume \thevolumenumber\ (\thevolumeyear)
\startpage--\finishpage\nl
Published: \publishdate}

\vglue .25truein

{\parskip=0pt\leftskip 0pt plus
1fil\def\\{\par\smallskip}{\Large\bf\thetitle}\par\medskip} \vglue
0.05truein

{\parskip=0pt\leftskip 0pt plus 1fil\def\\{\par}{\sc\theauthors}
\par\medskip}%
 
\vglue 0.03truein 

{\small\leftskip 25truept\rightskip 25truept{\bf Abstract}\stdspace\theabstract

{\bf AMS Classification}\stdspace\theprimaryclass
\ifx\thesecondaryclass\relax\else; \thesecondaryclass\fi\par
{\bf Keywords}\stdspace \thekeywords\par}\vglue 7truept

}   

\ifplaintex
\font\phead=cmsl9 scaled 950
\font\pnum=cmbx10 scaled 913
\font\pfoot=cmsl9 scaled 950
\headline{\vbox to 0pt{\vskip -4.5mm\line{\small\phead\ifnum
\count0=\startpage ISSN 1472-2739 (on-line) 1472-2747 (printed)
\hfill {\pnum\folio}\else\ifodd\count0\def\\{ }%
\ifx\theshorttitle\relax\thetitle\else\theshorttitle\fi\hfill{\pnum\folio}
\else\def\\{ and }{\pnum\folio}\hfill\ifx\theshortauthors\relax\theauthors
\else\theshortauthors\fi\fi\fi}\vss}}
\footline{\vbox to 0pt{\vglue 0mm\line{\small\pfoot\ifnum\count0=\startpage
\copyright\ \gtp\hfill\else
\agt, Volume \thevolumenumber\ (\thevolumeyear)\hfill\fi}\vss}}
\else
\font=cmsl9 scaled 1050
\font\lnum=cmbx10 
\font=cmsl9 scaled 1050
\makeatletter
\def\@oddhead{{\small\lhead\ifnum\count0=\startpage ISSN 1472-2739 
(on-line) 1472-2747 (printed)\hfill {\lnum\number\count0}\else\ifodd\count0
\def\\{ }\ifx\theshorttitle\relax \thetitle \else\theshorttitle\fi\hfill
{\lnum\number\count0}\else\def\\{ and }{\lnum\number\count0}
\hfill\ifx\theshortauthors\relax 
\theauthors\else\theshortauthors\fi\fi\fi}}\def\@evenhead{\@oddhead}
\def\@oddfoot{\small\lfoot\ifnum\count0=\startpage\copyright\ \gtp\hfill\else
\agt, Volume \thevolumenumber\ (\thevolumeyear)\hfill\fi}
\def\@evenfoot{\@oddfoot}
\makeatother
\fi
\let\maketitlepage\makeagttitle
\let\makeshorttitle\maketitlepage
\let\maketitle\maketitlepage

\newwrite\gtoutfile
\long\gdef\makeheadfile{  
{\def\\{, }\def\s{ }
\immediate\openout\gtoutfile head.xxx
\immediate\write\gtoutfile{To: math@arxiv.org}
\immediate\write\gtoutfile{Subject: put OR rep NNNNN:ppppp}
\immediate\write\gtoutfile{--text follows this line--}
\immediate\write\gtoutfile{Proxy-for: \ifx\theasciiauthors\relax
\theauthors\else\theasciiauthors\fi\s<\ifx\theasciiemail\relax\theemail\else\theasciiemail\fi>}
\immediate\write\gtoutfile{\noexpand\\}
\immediate\write\gtoutfile{Authors: \ifx\theasciiauthors\relax
\theauthors\else\theasciiauthors\fi}
{\def\\{ }\immediate\write\gtoutfile{Title: \ifx\theasciititle\relax
\thetitle\else\theasciititle\fi}}
\immediate\write\gtoutfile{Subj-class: GT or SG, GR etc}
\immediate\write\gtoutfile{MSC-class: \theprimaryclass\ifx\thesecondaryclass\relax\else, \thesecondaryclass\fi}
\immediate\write\gtoutfile{Journal-ref: Algebr. Geom. Topol. \thevolumenumber\s
(\thevolumeyear) \startpage-\finishpage}
\immediate\write\gtoutfile{Comments: Published by Algebraic and
Geometric Topology at}
\immediate\write\gtoutfile{\s\s\s  http://www.maths.warwick.ac.uk/agt/AGTVol\thevolumenumber/agt-\thevolumenumber-\thepapernumber.abs.html}
\immediate\write\gtoutfile{\noexpand\\}
\immediate\write\gtoutfile{}
\ifx\theasciiabstract\relax
\immediate\write\gtoutfile{\theabstract}\else
\immediate\write\gtoutfile{\theasciiabstract}\fi
\immediate\write\gtoutfile{}
\immediate\write\gtoutfile{\noexpand\\}
\immediate\write\gtoutfile{}
\immediate\closeout\gtoutfile}}  

\def\maketitlepage{\makeagttitle\makeheadfile}
\let\makeshorttitle\maketitlepage
\let\maketitle\maketitlepage

\def\ifplaintex{\expandafter\ifx\csname documentclass\endcsname\relax}

\def\gtp{{\mathsurround=0pt\it $\cal G\mskip-2mu$eometry \&\ 
$\cal T\!\!$opology $\cal P\!$ublications}}  

\def\recd{{\small Received:\qua\receiveddate\ifx\reviseddate\relax
\else\qquad Revised:\qua\reviseddate\fi\par}} 

\def\lognumber#1{\def\thelognumber{#1}}
\def\volumenumber#1{\def\thevolumenumber{#1}}
\def\volumeyear#1{\def\thevolumeyear{#1}}
\def\papernumber#1{\def\thepapernumber{#1}}
\def\pagenumbers#1#2{\def\startpage{#1}\def\finishpage{#2}}
\def\published#1{\def\publishdate{#1}}

\def\received#1{\def\receiveddate{#1}}
\def\revised#1{\def\reviseddate{#1}}
\def\accepted#1{\def\accepteddate{#1}}

\long\def\asciiabstract#1{\long\def\theasciiabstract{#1}}

\let\\\par\let\thelognumber\relax\let\thevolumenumber\relax
\let\thepapernumber\relax\let\thevolumeyear\relax\let\startpage\relax
\let\finishpage\relax\let\publishdate\relax\let\receiveddate\relax
\let\reviseddate\relax\let\accepteddate\relax\let\theasciititle\relax
\let\theasciiauthors\relax
\let\theasciiabstract\relax

\let\theasciiemail\relax

\ifplaintex
\font\logobig=cmssbx10 scaled 3836
\font\logomed=cmssbx10 scaled 2557
\else
\font\logobig=cmssbx10 scaled 4200
\font\logomed=cmssbx10 scaled 2800
\fi

\long\def\makeagttitle{   
\count0=\startpage
\agt\hfill      
\hbox to 45truept{\vbox to 0pt{\vglue -13truept{\logomed A\kern -.37em{\logobig 
T}\kern -.38em G}\vss}\hss}
\break
{\small Volume \thevolumenumber\ (\thevolumeyear)
\startpage--\finishpage\nl
Published: \publishdate}

\vglue .25truein

{\parskip=0pt\leftskip 0pt plus
1fil\def\\{\par\smallskip}{\Large\bf\thetitle}\par\medskip} \vglue
0.05truein

{\parskip=0pt\leftskip 0pt plus 1fil\def\\{\par}{\sc\theauthors}
\par\medskip}%
 
\vglue 0.03truein 

{\small\leftskip 25truept\rightskip 25truept{\bf Abstract}\stdspace\theabstract

{\bf AMS Classification}\stdspace\theprimaryclass
\ifx\thesecondaryclass\relax\else; \thesecondaryclass\fi\par
{\bf Keywords}\stdspace \thekeywords\par}\vglue 7truept

}   

\ifplaintex
\font\phead=cmsl9 scaled 950
\font\pnum=cmbx10 scaled 913
\font\pfoot=cmsl9 scaled 950
\headline{\vbox to 0pt{\vskip -4.5mm\line{\small\phead\ifnum
\count0=\startpage ISSN 1472-2739 (on-line) 1472-2747 (printed)
\hfill {\pnum\folio}\else\ifodd\count0\def\\{ }%
\ifx\theshorttitle\relax\thetitle\else\theshorttitle\fi\hfill{\pnum\folio}
\else\def\\{ and }{\pnum\folio}\hfill\ifx\theshortauthors\relax\theauthors
\else\theshortauthors\fi\fi\fi}\vss}}
\footline{\vbox to 0pt{\vglue 0mm\line{\small\pfoot\ifnum\count0=\startpage
\copyright\ \gtp\hfill\else
\agt, Volume \thevolumenumber\ (\thevolumeyear)\hfill\fi}\vss}}
\else
\font=cmsl9 scaled 1050
\font\lnum=cmbx10 
\font=cmsl9 scaled 1050
\makeatletter
\def\@oddhead{{\small\lhead\ifnum\count0=\startpage ISSN 1472-2739 
(on-line) 1472-2747 (printed)\hfill {\lnum\number\count0}\else\ifodd\count0
\def\\{ }\ifx\theshorttitle\relax \thetitle \else\theshorttitle\fi\hfill
{\lnum\number\count0}\else\def\\{ and }{\lnum\number\count0}
\hfill\ifx\theshortauthors\relax 
\theauthors\else\theshortauthors\fi\fi\fi}}\def\@evenhead{\@oddhead}
\def\@oddfoot{\small\lfoot\ifnum\count0=\startpage\copyright\ \gtp\hfill\else
\agt, Volume \thevolumenumber\ (\thevolumeyear)\hfill\fi}
\def\@evenfoot{\@oddfoot}
\makeatother
\fi
\let\maketitlepage\makeagttitle
\let\makeshorttitle\maketitlepage
\let\maketitle\maketitlepage

\newwrite\gtoutfile
\long\gdef\makeheadfile{  
{\def\\{, }\def\s{ }
\immediate\openout\gtoutfile head.xxx
\immediate\write\gtoutfile{To: math@arxiv.org}
\immediate\write\gtoutfile{Subject: put OR rep NNNNN:ppppp}
\immediate\write\gtoutfile{--text follows this line--}
\immediate\write\gtoutfile{Proxy-for: \ifx\theasciiauthors\relax
\theauthors\else\theasciiauthors\fi\s<\ifx\theasciiemail\relax\theemail\else\theasciiemail\fi>}
\immediate\write\gtoutfile{\noexpand\\}
\immediate\write\gtoutfile{Authors: \ifx\theasciiauthors\relax
\theauthors\else\theasciiauthors\fi}
{\def\\{ }\immediate\write\gtoutfile{Title: \ifx\theasciititle\relax
\thetitle\else\theasciititle\fi}}
\immediate\write\gtoutfile{Subj-class: GT or SG, GR etc}
\immediate\write\gtoutfile{MSC-class: \theprimaryclass\ifx\thesecondaryclass\relax\else, \thesecondaryclass\fi}
\immediate\write\gtoutfile{Journal-ref: Algebr. Geom. Topol. \thevolumenumber\s
(\thevolumeyear) \startpage-\finishpage}
\immediate\write\gtoutfile{Comments: Published by Algebraic and
Geometric Topology at}
\immediate\write\gtoutfile{\s\s\s  http://www.maths.warwick.ac.uk/agt/AGTVol\thevolumenumber/agt-\thevolumenumber-\thepapernumber.abs.html}
\immediate\write\gtoutfile{\noexpand\\}
\immediate\write\gtoutfile{}
\ifx\theasciiabstract\relax
\immediate\write\gtoutfile{\theabstract}\else
\immediate\write\gtoutfile{\theasciiabstract}\fi
\immediate\write\gtoutfile{}
\immediate\write\gtoutfile{\noexpand\\}
\immediate\write\gtoutfile{}
\immediate\closeout\gtoutfile}}  

\def\maketitlepage{\makeagttitle\makeheadfile}
\let\makeshorttitle\maketitlepage
\let\maketitle\maketitlepage

\lognumber{23}
\volumenumber{2}
\volumeyear{2002}
\papernumber{23}
\published{22 June 2002}
\pagenumbers{465}{497}
\received{24 October 2001}
\revised{24 May 2002}
\accepted{6 June 2002}

\usepackage{graphicx,amssymb,subfigure}

\def\lacute{\mathopen{<}}
\def\racute{\mathclose{>}}
\def\defeq{ \stackrel{\mbox{{\tiny def}}}{=} }
\def\finpf{\endproof}
\def\defm{{\mathcal{H}(M)}}
\def\param{{{\cal P}(K)}}

\def\isom{{\mbox{PSL}_2\mathbb{C}}}
\def\Char{{X(M)}}
\def\tM{{ \widetilde{M} }}
\def\tN{{ \widetilde{N} }}
\def\dist{{{\mit d}_{\mathbb{C}}}}
\def\hyp{{ \mathbb{H}^3 }}
\def\tr{{ \mbox{tr} }}
\def\si{{\mathbb{S}^2_\infty}}

\def\orth{{ \mathit{Orth}_K }}
\def\hol{{ \mathit{Hol} }}
\def\orthz{{ \mathit{Orth}_{{K_0}} }}
\def\OrthF{{ D_F }}
\def\Od{{ D_\triangle }}
\def\Pd{{{\cal P}(\triangle)}}
\def\rk{{ rank }}
\def\transpose{{ ^{\mbox{{\tiny T}}} }}
\def\endpts{{ \si \times \si \setminus \Delta }}
\def\OrLines{{ \mbox{S} \mathcal{L} }}
\def\Lines{{ \mathcal{L} }}
\def\rep{{\mathcal{R}(M)}}
\def\redh{{\mathcal{RH}_4}}
\def\rhex{{\mathcal{RH}_3}}

\def\dom{{\mathcal{D}}}
\def\mbar{{ \overline{M} }}

\newtheorem{theorem}{Theorem}
\newtheorem{prop}[theorem]{Proposition}
\newtheorem{cor}[theorem]{Corollary}
\newtheorem{lemma}[theorem]{Lemma}
\newtheorem{conj}[theorem]{Conjecture}

\theoremstyle{definition}
\newtheorem{defn}[theorem]{Definition}
\newtheorem{problem}[theorem]{Problem}

\begin{document}
\title{A new invariant on hyperbolic Dehn surgery space}
\author{James G. Dowty}
\address{Department of Mathematics,
University of Melbourne\\Parkville, 3052,
Australia}
\email{jamesdowty@bigpond.com.au}

\begin{abstract}
In this paper we define a new invariant of the incomplete hyperbolic
structures on a $1$-cusped finite volume hyperbolic $3$-manifold $M$,
called the ortholength invariant.  We show that away from a (possibly
empty) subvariety of excluded values this invariant both locally
parameterises equivalence classes of hyperbolic structures and is a
complete invariant of the Dehn fillings of $M$ which admit a
hyperbolic structure.  We also give an explicit formula for the
ortholength invariant in terms of the traces of the holonomies of
certain loops in $M$.  Conjecturally this new invariant is intimately
related to the boundary of the hyperbolic Dehn surgery space of
$M$.
\end{abstract}

\asciiabstract{%
In this paper we define a new invariant of the incomplete hyperbolic
structures on a 1-cusped finite volume hyperbolic 3-manifold M, called
the ortholength invariant.  We show that away from a (possibly empty)
subvariety of excluded values this invariant both locally
parameterises equivalence classes of hyperbolic structures and is a
complete invariant of the Dehn fillings of M which admit a hyperbolic
structure.  We also give an explicit formula for the ortholength
invariant in terms of the traces of the holonomies of certain loops in
M.  Conjecturally this new invariant is intimately related to the
boundary of the hyperbolic Dehn surgery space of M.}


\maketitle

\section{Introduction}
\label{S:Intro}

Let $M$ be an orientable 3-manifold which admits a complete, finite-volume
hyperbolic structure with a single cusp.
Thurston's hyperbolic Dehn surgery theorem applied to $M$ says that
all but a finite number of topological
Dehn fillings on $M$ have hyperbolic structures.
To prove this, Thurston introduced the deformation space $\defm$ of
incomplete hyperbolic structures on $M$ whose completions have `Dehn 
surgery-type' singularities (see \cite{ThurstonNotes}).   This space includes
any hyperbolic cone-manifold whose set of non-singular points is 
diffeomorphic to $M$. The space $\defm$ is of continuing interest because it
offers a possible approach (due to Thurston) to the Geometrization Conjecture.
In particular, non-hyperbolic geometric structures and topological
decompositions along incompressible spheres and tori can often be
produced by understanding
the ways that hyperbolic structures can degenerate near the boundary of
$\defm$ (see Thurston \cite[Chapter 4]{ThurstonNotes},
Kerckhoff \cite{Kerckhoff} and Kojima \cite{KojSurvey}).
This approach has recently yielded a proof of the Orbifold
Theorem (see Cooper-Hodgson-Kerckhoff \cite[Chapter 7]{OrbifoldNotes} or 
Boileau-Porti \cite{BPorti}).

In this paper we introduce a $\mathbb{C}^n$-valued invariant $\orth$ of each 
incomplete hyperbolic structure in $\defm$, called the {\em ortholength invariant}.
This invariant is defined in terms of a topological ideal 
triangulation $K$ of $M$, but it is actually independent
of $K$ in the sense that for `generic' $K$, $\orth$
determines the ortholength invariant corresponding to any other
ideal triangulation (see Theorem \ref{T:Param}).  The main result of this paper says that away from a (possibly empty) subvariety of values, 
the ortholength invariant locally parameterises $\defm$ and it is 
a complete invariant of the (topological) Dehn fillings of $M$ which admit 
a hyperbolic structure (see Theorem \ref{T:CP}). 

For a hyperbolic structure whose metric completion is a cone-manifold, the
ortholength invariant is essentially a list of hyperbolic cosines of the
complex  distances\footnote{The complex distance between two lines is the
hyperbolic distance between them plus $i$ times an angle of twist.} from the
cone-manifold's singular set to itself along the edges\footnote{Not all 
homotopy classes of paths from the singular set to itself contain a 
distance-realising geodesic, 
however the ortholength invariant is defined purely in terms of the holonomy 
representation and so is well-defined in any case.  } of $K$.  
Hence there is a close connection between the ortholength invariant and the 
{\em tube radius} of a hyperbolic Dehn filling, i.e.\ the  supremum of 
the radii of the embedded hyperbolic tubes in $M$ about the singular set of 
the Dehn filling.

The importance of the tube radius has emerged recently from the work of
Kojima \cite{KojPi} and Hodgson-Kerckhoff \cite{Kerckhoff}, \cite{HK}.
In these papers,
the condition that the tube radius stays bounded away from zero as cone-angles
are varied has emerged as the key to ensuring that the hyperbolic
structures do not degenerate.   (By contrast, it is possible for
hyperbolic cone-manifolds to degenerate while their volume remains bounded 
above zero, e.g.\ see \cite{KojPi}, especially Theorem 7.1.2 or the 
example of \S 7.2.)
This suggests a relationship between the boundary of $\defm$  and
the ortholength invariant which deserves further study (see Conjecture 
\ref{C:Conjecture}).  

The terminology `ortholength invariant' follows  \cite{Meyerhoff}, where 
Meyerhoff defines the complex ortholength 
spectrum of a closed hyperbolic 3-manifold as the set of complex distances 
between pairs of the manifold's simple closed
geodesics.  He shows that the ortholength 
spectrum plus some combinatorial data determines the manifold up to isometry. 
Hence the fact that $\orth$ is a complete invariant of the topological 
Dehn fillings of $M$ which admit a hyperbolic structure 
can be interpreted as a strengthened version 
of Meyerhoff's result, namely that a finite (and computable) subset of the
ortholength spectrum plus slightly stronger 
combinatorial information than assumed in \cite{Meyerhoff} determines 
the filled manifold up to isometry.  

The rest of this paper is set out as follows.  In Section \ref{S:Lines} we 
show that it is possible to parameterise a configuration
of lines in $\hyp$ in terms of the complex distances  between them.  
In Section \ref{S:OrthInv} we define the ortholength invariant 
and give a formula for it in terms of the traces of the holonomies of certain
loops in $M$, thereby showing that $\orth$ is a rational map from the 
$\isom$-character variety $\Char$ of $M$ into $\mathbb{C}^n$.  
The main section of this paper is Section \ref{S:Param}, where we prove 
that the ortholength invariant locally parameterises incomplete 
hyperbolic structures on $M$ and is a complete invariant when restricted 
to the (topological) Dehn fillings of $M$ which admit a hyperbolic structure 
(for those hyperbolic structures whose ortholength invariant does not lie in a 
certain subvariety).  We also prove  
(under a weak technical assumption) that the ortholength invariant 
is a birational equivalence from certain interesting irreducible 
components of $\Char$ to certain irreducible components of a complex 
affine algebraic variety $\param \subseteq \mathbb{C}^n$.  
In Section \ref{S:Eg} we calculate the ortholength invariant when $M$ is 
the figure-8 knot complement.  
Section \ref{S:Conc} concludes this paper with a conjecture and a discussion 
of some applications for the ortholength invariant, including 
the construction of incomplete hyperbolic structures on $M$.

This paper is based on my doctoral thesis and I would like to acknowledge 
and thank my advisor Craig D. Hodgson for his guidance throughout the work presented here.

\section{Configurations of lines in hyperbolic 3-space}
\label{S:Lines}

We say that $n$ oriented lines  $l_1, \ldots, l_n$ in
hyperbolic 3-space $\hyp$ {\em realise} a set of complex numbers 
$x_{ij} \in \mathbb{C}$ ($i,j = 1, \ldots ,n$) if 
$\cosh$ of the complex
distance\footnote{The complex distance between two lines is the hyperbolic
distance between them plus $i$ times an angle of twist---see Fenchel 
\cite[\S V.3]{Fenchel}.}
between each pair of lines $l_i$ and $l_j$ is equal to $x_{ij}$.
This section is motivated by the following question.  

\begin{problem}
\label{Q:InfCh1}
When is it possible to realise a given set of complex numbers 
$x_{ij} \in \mathbb{C}$  ($i,j = 1, \ldots ,n$) 
by an arrangement of lines and how unique is such an arrangement?  
\end{problem}
A complete answer to Problem \ref{Q:InfCh1} is given in 
Theorem \ref{T:SolnCh1}.  A corollary to this theorem (Corollary \ref{C:N=4}) 
says that there exists
an arrangement of four lines which realises a set of complex numbers 
$x_{ij} \in \mathbb{C}$ ($i,j = 1, \ldots ,4$) if 
and only if these complex numbers satisfy a `hextet' equation.
This corollary is the key result which allows us to locally parameterise 
hyperbolic structures on a $3$-manifold by `ortholengths' 
(see Section \ref{S:Param}).  
The results of this section prior to Problem \ref{Prob:Ch1} are simply 
re-statements of some of the ideas Fenchel 
presented in \cite{Fenchel}.   For related material see
Thurston \cite[\S\S 2.3-2.6]{ThurstonBook}.

The end-points of a line in $\hyp$ are two distinct points on the sphere 
at infinity $\si$, and conversely any two such points determine a line.  
We identify the set of {\em oriented lines} in $\hyp$ with  
$\endpts$ where $\Delta = \{ (x,x) \mid x \in \si\}$, and we refer to 
an element of $\endpts$  as the {\em ordered end-points} of an oriented 
line in $\hyp$.  

A simple calculation shows that every traceless matrix
$A \in \mbox{SL}_2\mathbb{C}$ has eigenvalues $\pm i$.  The eigenvectors of 
$A$ correspond to points on the sphere at infinity $\si$ which
are fixed by the action\footnote{%
$A$ acts as a projective transformation on $\si = \mathbb{CP}^1$.} of $A$.  
Hence there is an oriented line corresponding to $A$ whose ordered 
end-points are $(p,q) \in \endpts$, where $p$ corresponds to eigenvalue $-i$ 
and $q$ corresponds to eigenvalue $i$.  It is not hard to check that 
the map $A \mapsto (p,q)$ is a homeomorphism from 
$$\OrLines \defeq \{ l \in \mbox{SL}_2\mathbb{C} \mid \mbox{tr}~ l = 0 \}$$
to the space of oriented lines in $\hyp$.  
We identify these two spaces and from now on we will use the same symbol to
denote both an oriented line and the corresponding element of $\OrLines$.

Note that as a consequence of this, $-l \in \OrLines$ denotes the same 
line in $\hyp$ as $l \in \OrLines$ but with the opposite orientation.

Now, $\isom$ acts by orientation-preserving isometries on the set 
of oriented lines in $\mathbb{H}^3$.
Under the correspondence between oriented lines and $\OrLines$, this  
gives us a corresponding action of $\isom$ on $\OrLines$, given 
by $g \cdot l = \tilde{g} l \tilde{g}^{-1}$ for any $g \in \isom$ and 
$l \in \OrLines$, where $\tilde{g} \in 
\mbox{SL}_2\mathbb{C}$ is either of the matrices covering $g$.  
Hence studying the geometric properties of arrangements of lines in 
$\hyp$ is the 
same as studying the properties of subsets of $\OrLines$ which are invariant 
under the conjugacy action of $\mbox{SL}_2\mathbb{C}$. 

Given two oriented lines $l,m \in \OrLines$, an obvious conjugacy invariant
of the pair is $\mbox{tr} (lm)$.  Lemma \ref{L:InnerProduct} (below) says that 
this invariant is essentially $\cosh$ of the complex distance between 
$l$ and $m$.  However, before stating this lemma we introduce some 
notation\footnote{\label{FootnoteRef}I am grateful to the referee for drawing  
my attention to the fact that $\Lines$ is the Lie algebra of 
$\mbox{SL}_2\mathbb{C}$, the corresponding action of $\mbox{SL}_2\mathbb{C}$ on 
$\Lines$ is the adjoint action, 
the form $\lacute \cdot,\cdot \racute$ (below) is a 
multiple of the Killing form and Lemma \ref{L:IsomL} (below) is essentially 
the (known) result that $\isom$ is isomorphic to $\mbox{SO}_3\mathbb{C}$.
}. 

\begin{defn}[Line Matrices]
The complex vector space $\Lines$ of {\em line matrices} consists of
all $2 \times 2$ complex matrices with zero trace.  This space is endowed with 
a symmetric bilinear form $\lacute \cdot,\cdot \racute$ given by
$$ \lacute l,m \racute = - \frac{1}{2} \mbox{tr} (lm)$$
for any $l,m \in \Lines$. 
\end{defn}

It is not hard to check that the bilinear form $\lacute \cdot,\cdot \racute$ 
is non-degenerate on $\Lines$.  
Also, each non-singular line matrix $l \in \Lines$ acts on 
$\mathbb{H}^3$ as a half-turn about some line so there is an 
(unoriented) line associated with each 
non-singular line matrix.  This is where the terminology `line matrices' 
comes from.  

\begin{lemma}
\label{L:InnerProduct}
For any $l, m \in \OrLines$, 
\begin{eqnarray}
\label{E:FundLink}
\lacute l,m \racute = \cosh ( \dist (l,m)) 
\end{eqnarray}
where $\dist (l,m)$ is the complex distance between $l$ and $m$. 
\end{lemma}

\proof See Fenchel \cite[\S V.3]{Fenchel}.
\finpf

Alternatively, the reader can take Lemma \ref{L:InnerProduct} as the 
definition of the complex distance between two oriented lines in $\hyp$.

Note that $\lacute l,l \racute = \det l$ for any $l \in \Lines$
so $\OrLines = \{l\in \Lines \mid \lacute l,l \racute = 1 \}$
is the set of normalised line matrices.
Then by Lemma \ref{L:InnerProduct} we can rephrase Problem \ref{Q:InfCh1} 
in terms of linear algebra as follows.  

\begin{problem}
\label{Prob:Ch1}
Given some complex numbers $x_{ij} \in \mathbb{C}$  ($i, j = 1, \ldots , n$) 
so that $x_{ij} = x_{ji}$ and $x_{ii} = 1$, do there exist oriented lines 
$l_1, \ldots, l_n \in \OrLines$ so that 
$$\lacute l_i, l_j \racute = x_{ij}?$$
To what extent do these conditions determine the lines $l_1, \ldots, l_n$ 
if they exist?
\end{problem}

Now, the action of $\isom$ on $\OrLines$ extends to an action on 
$\Lines$ in an obvious way\footnote{For $g \in \isom$ and $l\in\Lines$, 
$g\cdot l =
\tilde{g} l \tilde{g}^{-1}$ where $\tilde{g} \in \mbox{SL}_2\mathbb{C}$
is one of the two matrices covering $g$.} and clearly $\lacute \cdot,\cdot 
\racute$ is invariant under this action.   Also note that 
$\lacute \cdot,\cdot \racute$ is invariant under the map 
$\Lines \to \Lines$ given by $l \mapsto -l$.  
The following lemma says that these are the only isomorphisms of $(\Lines, 
\lacute \cdot,\cdot \racute)$.  

\begin{lemma}
\label{L:IsomL}
For each linear map $\phi: \Lines \to \Lines$ which preserves 
$\lacute \cdot,\cdot \racute$ there is some element $g \in \isom$ and some 
choice $\pm 1$ of sign so that 
$$ \pm \phi (l)  = g \cdot l $$
for any $l \in \Lines$.  
\end{lemma}

\begin{figure}[ht!]
\vspace{-0.5cm}
\begin{center}
\subfigure[Left-handed]{
\includegraphics[height=4cm]{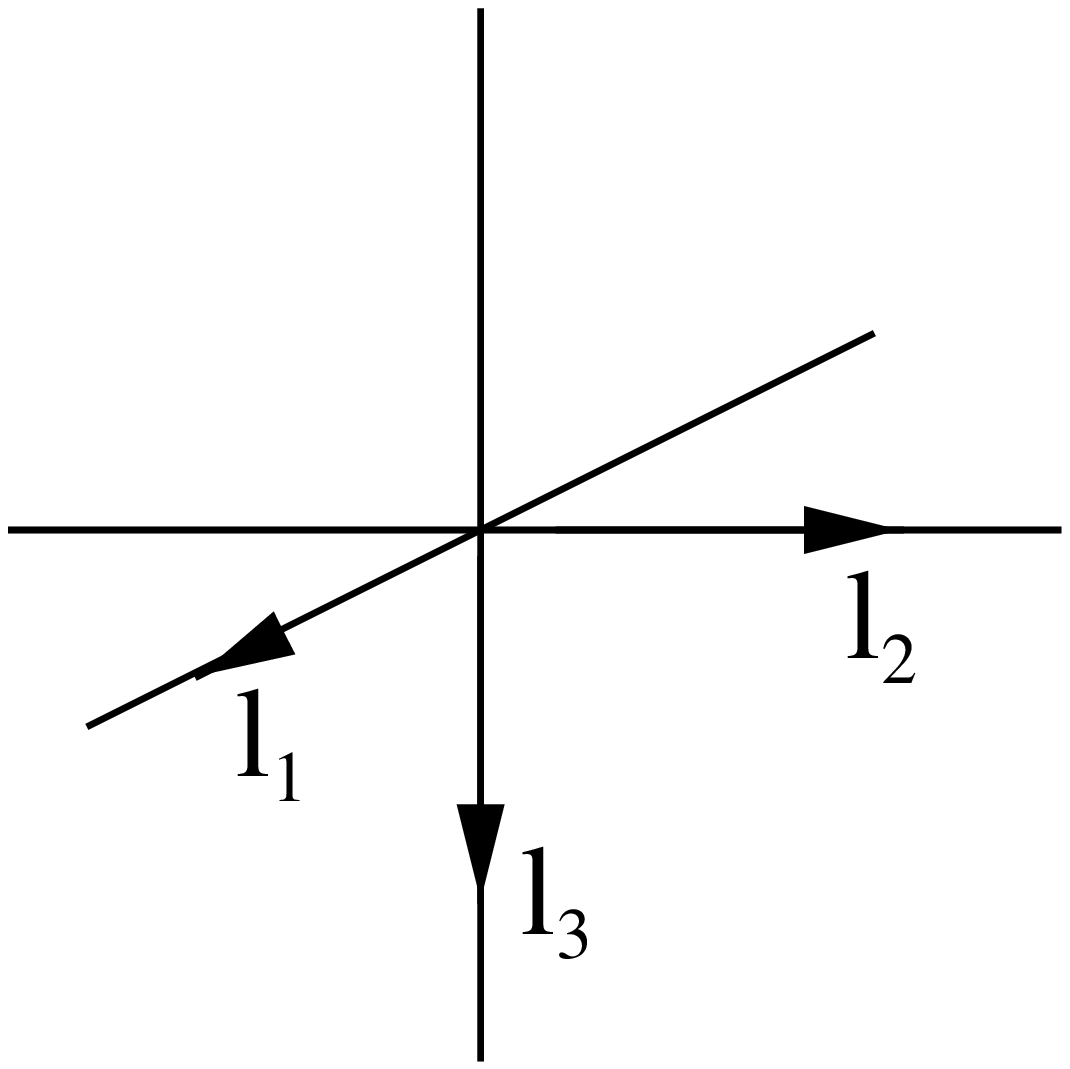}\hspace{1cm}
}
\subfigure[Right-handed]{
\includegraphics[height=4cm]{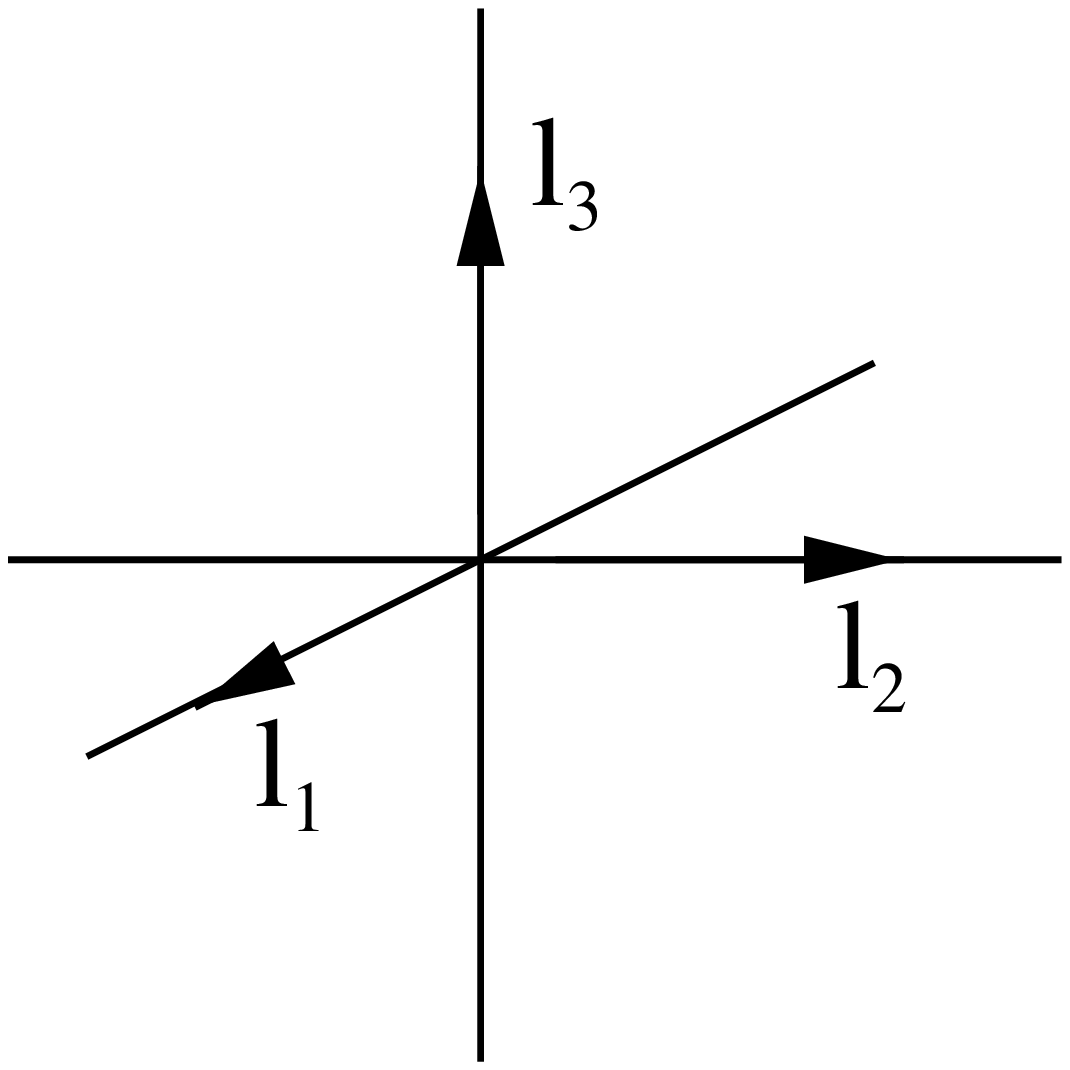}
} 
\end{center}\vspace{-0.5cm}
\caption{The two orientation-classes of orthogonal frames} 
\label{F:CoordSystem}
\end{figure}

\proof
Consider three lines $l_1,l_2,l_3 \in \OrLines$ for which  
$\lacute l_i,l_j\racute = 0$ for each $i\not= j$.  
By Lemma \ref{L:InnerProduct}, 
these lines all  meet at a common point of $\hyp$ where they are mutually
perpendicular. We assume that the orientations of the $l_i$ have been chosen 
so that they define a right-handed frame at this common point 
(see Figure \ref{F:CoordSystem}).  

Now, if $\phi: \Lines \to \Lines$ preserves $\lacute \cdot,\cdot \racute$
then $\lacute \phi(l_i), \phi(l_j) \racute = 0$ for each $i\not= j$ and so 
by Lemma \ref{L:InnerProduct},
$\phi(l_1), \phi(l_2), \phi(l_3)$ are also mutually perpendicular lines 
which meet at a common point.  
By replacing $\phi$ by $- \phi$ if need be we can assume
$\phi(l_1), \phi(l_2), \phi(l_3)$ define a right-handed frame. 
Since $\isom$ acts transitively on the bundle of
right-handed orthogonal frames of $\hyp$ (see \cite[\S 2.2]{ThurstonBook})
there is an isometry $g\in \isom$ which takes $l_i$ to $\phi(l_i)$ for each 
$i = 1,2,3$.  Hence $\phi(l_i)  = g \cdot l_i$ for each $i = 1,2,3$. 
Since the $l_i$ form a basis for $\Lines$ (linear dependence would imply that 
$\lacute l_i,l_i\racute = 0$ for some $i$) this proves the lemma. 
\finpf

In fact, from this proof it is clear that the action of $\isom \times \{\pm1\}$
on $\Lines$ induces an isomorphism from $\isom \times \{\pm1\}$ to 
the group $\mbox{O}_3(\mathbb{C})$ of isomorphisms of 
$(\Lines,\lacute \cdot,\cdot \racute)$ (see also Footnote \ref{FootnoteRef}, 
above).  Note that $-1$ acts on $\Lines$ by simultaneously reversing the 
orientations of all lines 
and is not related to the orientation-reversing isometries of $\hyp$.  

Let $(V,g)$ denote a finite-dimensional complex vector space $V$ equipped 
with a symmetric, bilinear form $g$.   Then associated to $g$ there is a map 
$V \to V^*$ from $V$ to its dual  given by $v \mapsto g(v,\cdot)$.  
The {\em rank} of $g$ (denoted $\rk(g)$) is the dimension of 
the image of this map and $g$ is {\em non-degenerate} if this 
map is an isomorphism.  To solve Problem \ref{Prob:Ch1} we will use the 
following lemma from linear algebra.  

\begin{lemma}
\label{L:BiForm}
Let $(V,g)$ and $(W,h)$ be two finite-dimensional 
complex vector spaces equipped with symmetric, bilinear forms $g$ and $h$.
If $\rk(g) \leq \rk(h)$ then there exists a linear transformation 
$\phi : V \to W$ so that
$$\phi^*h = g$$
i.e.\ so that $g(x,y) = h(\phi x, \phi y)$ for any $x,y \in V$.  Furthermore, 
if $h$ is non-degenerate and $\rk(g) = \rk(h)$ or $\rk(g) = \rk(h)-1$ 
then $\phi$ is unique up to composing it on the left with an isomorphism 
of $(W,h)$.   
\end{lemma}

\proof  This is an elementary consequence of the fact that 
$(V,g)$ is isomorphic to $(\mathbb{C}^n, E_r)$ for $n = \mbox{dim}(V)$ and 
$r = \rk(g)$, where $E_r$ is the 
bilinear form defined by $E_r(u,v) = u\transpose Gv$ for any $u,v \in 
\mathbb{C}^n$ and $G$ is the 
$n \times n$ diagonal matrix $\mbox{diag}(1, \ldots, 1, 0, \ldots, 0)$ 
with $r$ non-zero entries.  
\finpf

Armed with Lemmas \ref{L:IsomL} and \ref{L:BiForm} we can now completely solve
Problem \ref{Prob:Ch1}.  

\begin{theorem}
\label{T:SolnCh1}
For $i,j = 1,\ldots , n$ let $x_{ij}  \in \mathbb{C}$ be given
complex numbers  so that $x_{ij} = x_{ji}$ and $x_{ii} =1$.  
Then there exist $l_1, \ldots, l_n \in \OrLines$ for which 
$\lacute l_i, l_j \racute  = x_{ij}$ if and only if $\rk(X) \leq 3$, where 
$X$ is the $n \times n$ matrix whose $(i,j)^{th}$ entry is $x_{ij}$.  If some  
$x_{ij} \not= \pm 1$ then the arrangement of lines $l_1, \ldots, l_n$ 
is unique up to the action of $\isom \times \{\pm 1\}$, 
i.e.\ unique up to orientation-preserving isometry and simultaneous 
reversal of orientations.
\end{theorem}

If all $x_{ij} = \pm 1$ then there may be non-isometric arrangements of lines 
which realise the $x_{ij}$.

\proof[Proof of Theorem \ref{T:SolnCh1}]  
Let $x_{ij} \in \mathbb{C}$ be given complex numbers and
suppose that there exist lines $l_1, \ldots, l_n \in \OrLines$ so that
$x_{ij} =  \lacute l_i, l_j \racute$.
Then since $\Lines$ is 3-dimensional, any four of the lines must be linearly
dependent.   From this it follows that any four rows of $X$ are linearly 
dependent. For example,  there exist non-zero constants 
$\alpha_1, \ldots, \alpha_4 \in \mathbb{C}$  so that 
$\alpha_1 l_1 + \ldots + \alpha_4 l_4 = 0$ and hence  
$$\alpha_1 \lacute l_1, l_i \racute +\ldots+ \alpha_4 \lacute l_4, l_i \racute= 0$$ 
(i.e.\ $\alpha_1 x_{1i}+\ldots+\alpha_4 x_{4i} =0$) for each $i = 1, \ldots, n$.
Hence $\rk(X) \leq 3$.  

So now suppose that for $i,j = 1, \ldots, n$ we have 
$x_{ij} \in \mathbb{C}$ so that $x_{ij} = x_{ji}$ and $x_{ii}=1$ and 
that $\rk(X) \leq 3$ where $X = [x_{ij}]$ is the 
matrix defined in the statement. 
Define a symmetric, bilinear form $g$ on $\mathbb{C}^n$  by 
$$g(x, y) = x \transpose  X y $$
for any  $x, y \in \mathbb{C}^n$.  Then $\rk(g) \leq  3 = 
\rk (\lacute \cdot, \cdot \racute )$ so setting  
$(V,g) = (\mathbb{C}^n,g)$ and $(W,h) = (\Lines,\lacute \cdot,\cdot\racute)$ 
in Lemma \ref{L:BiForm} gives us a  
linear map $\phi : \mathbb{C}^n \to \Lines$ so that 
$\phi^* \lacute \cdot, \cdot \racute = g$.  
For each $i = 1, \ldots, n$ let $l_i \defeq \phi e_i$  where $e_1, \ldots, e_n$
is the standard basis for $\mathbb{C}^n$.   Then 
$$x_{ij} = e_i\transpose X e_j = g(e_i, e_j) =
\lacute \phi e_i, \phi e_j \racute = \lacute l_i, l_j \racute$$
for any $i,j = 1, \ldots, n$.  Note that 
$\lacute l_i, l_i \racute = x_{ii} = 1$ for each $i = 1, \ldots n$ 
so $l_i \in \OrLines$, i.e.\  
each $l_i$ corresponds to an oriented line in $\hyp$.  

Now, if some $x_{ij} \not= \pm 1$ then $\rk(g) = 2$ or $3$ (this uses 
the fact that $X$ is symmetric and $x_{ii} = 1$).
Hence by  Lemma \ref{L:BiForm}, $\phi$ (and hence the arrangement of lines 
$l_1, \ldots, l_n$) is unique up to composing $\phi$ 
by an isomorphism of $(\Lines,\lacute \cdot,\cdot\racute)$.  
But by Lemma \ref{L:IsomL} the isomorphisms of 
$(\Lines,\lacute \cdot,\cdot\racute)$ are exactly given by the 
action of $\isom \times \{\pm 1\}$ on $\Lines$.  
Hence the arrangement of lines $l_1, \ldots, l_n$ is unique up to 
orientation-preserving isometry and simultaneous reversal of orientations.  
\finpf

\begin{figure}[ht!]
\vspace{-0.5cm}
\begin{center}
\subfigure[A right-angled hexagon]{
\label{F:RightHexagon}
\includegraphics[height=4.5cm]{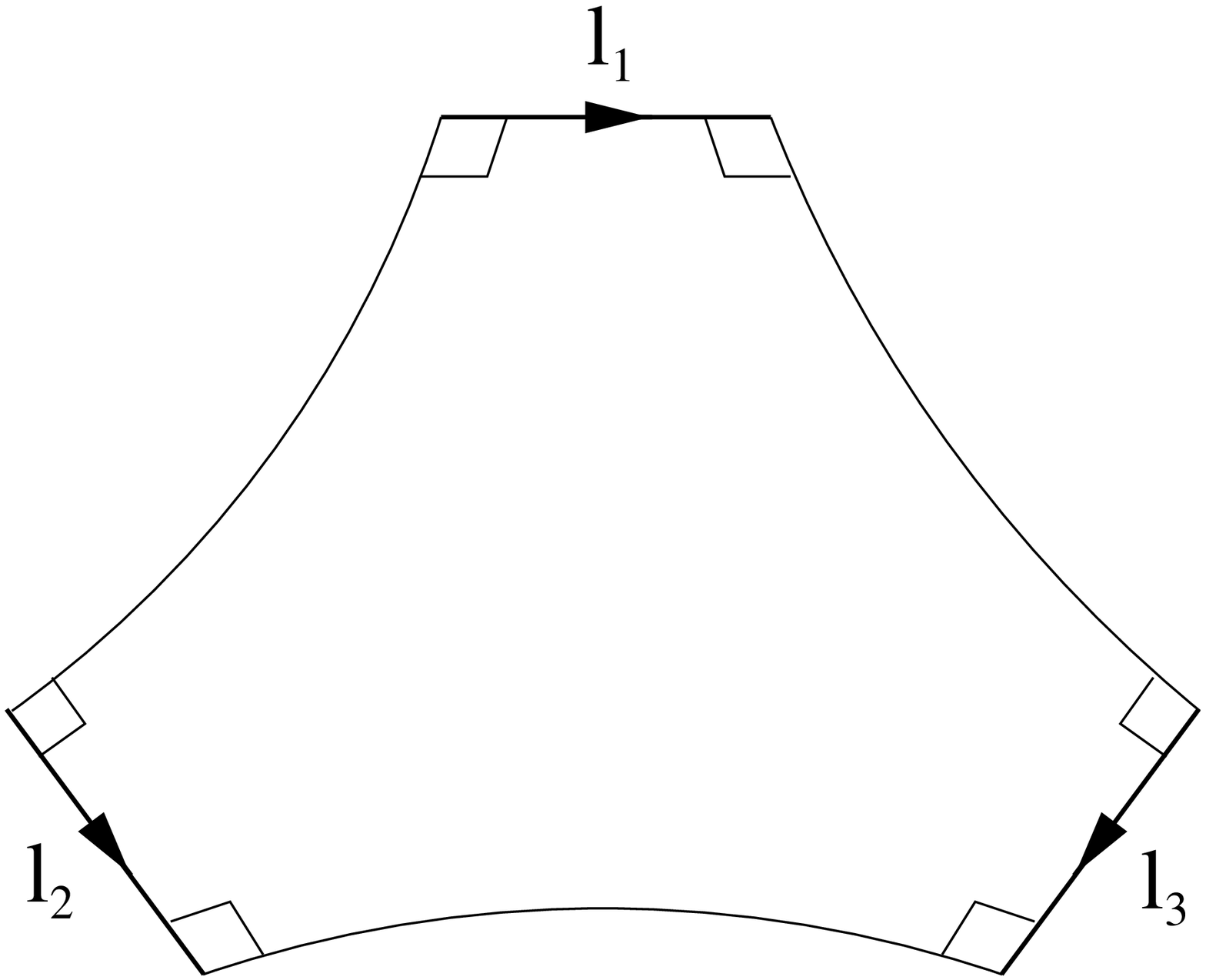}\hspace{1cm}
}
\subfigure[A hextet]{
\label{F:HexTet}
\includegraphics[height=4.5cm]{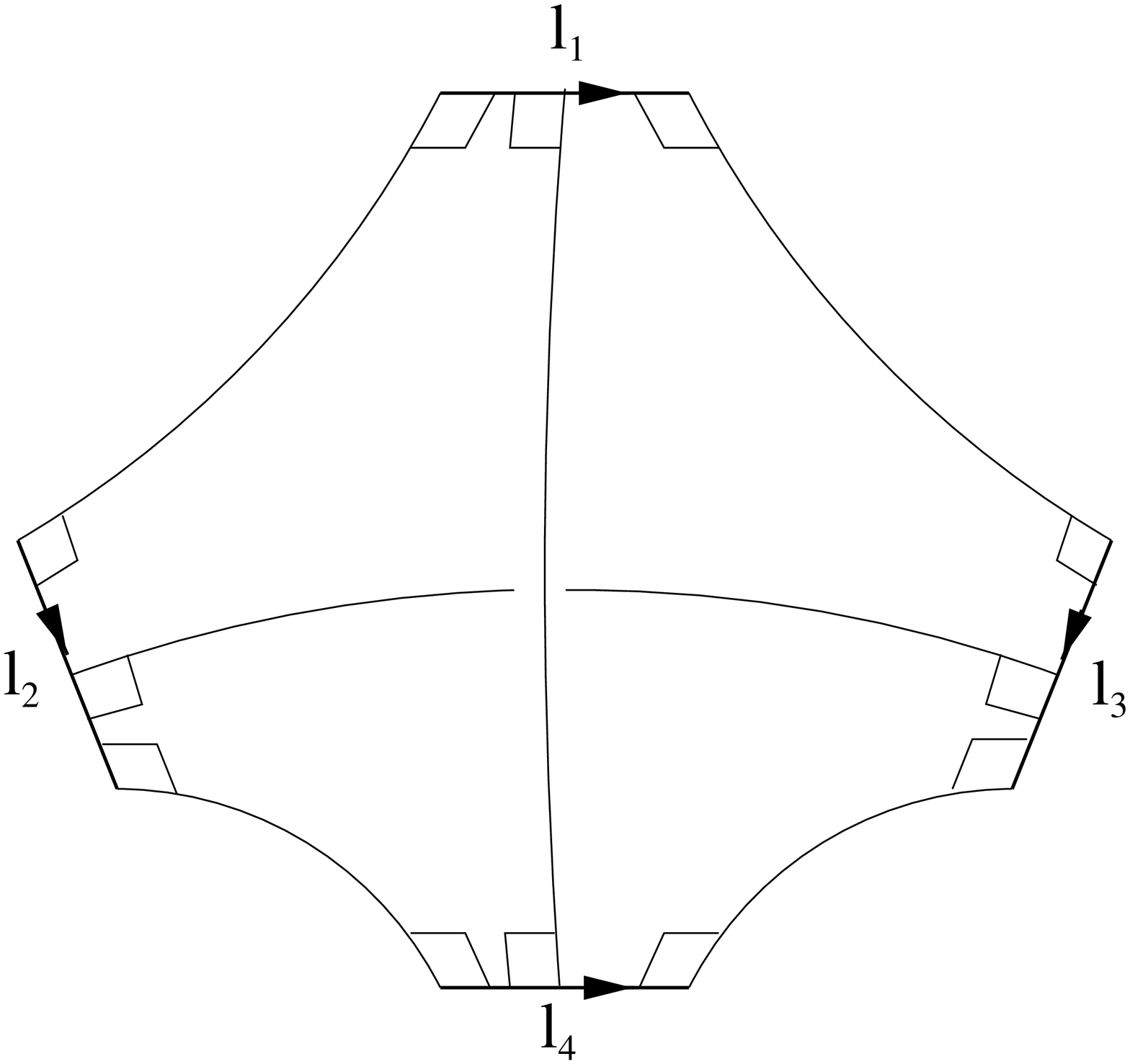}
} 
\end{center}\vspace{-0.5cm}
\caption{Arrangements of (a) three lines and (b) four lines}
\end{figure}

The special cases of Theorem \ref{T:SolnCh1} when $n= 3$ and $n=4$ are of 
particular interest to us.   If three lines $l_1,l_2,l_3 \in \OrLines$ all 
have distinct end-points on $\si$ then any two of these lines has a 
common-perpendicular.   Adding these three perpendiculars and truncating 
appropriately gives an arrangement 
of six line segments in $\hyp$ which meet at right-angles, i.e.\ a right-angled 
hexagon (see Figure \ref{F:RightHexagon}).  Motivated by this generic case
we will simply say that any three 
oriented lines define a right-angled hexagon, without making 
the assumption that the lines have distinct end-points.  Then Theorem 
\ref{T:SolnCh1} applied to $l_1,l_2,l_3$ gives the well-known result 
 that the complex distances along alternating edges 
of a right-angled hexagon determine the hexagon up to isometry (see 
\cite[\S VI.4]{Fenchel}).

Similarly, any four lines $l_1,\ldots,l_4 \in \OrLines$ which 
all have distinct end-points on $\si$ have six pair-wise 
common-perpendiculars.  Adding these 
perpendiculars and truncating all lines appropriately 
gives an arrangement of lines in $\hyp$ loosely resembling a 
tetrahedron (see Figure \ref{F:HexTet}).  
Since this arrangement is like a tetrahedron whose vertices 
have been stretched into the lines $l_1,\ldots,l_4$, turning its
faces into right-angled hexagons, we call 
such an arrangement a {\em hextet}.  As with right-angled hexagons, 
we will drop the 
requirement that the end-points of $l_1,\ldots,l_4$ be distinct and 
simply say that any four oriented lines in $\hyp$ define a hextet.  
Substituting $n=4$ into Theorem \ref{T:SolnCh1} gives the following
corollary.  

\begin{cor}[Existence and Rigidity of Hextets]
\label{C:N=4}
Let $x_{ij}  \in \mathbb{C}$ be given
complex numbers for $i,j = 1,\ldots , 4$ so that $x_{ij} = x_{ji}$, $x_{ii} =1$.  
Then there exist $l_1, \ldots, l_4 \in \OrLines$ so that
$\lacute l_i, l_j \racute  = x_{ij}$ if and only if the $x_{ij}$ satisfy the
{\em hextet equation}
\begin{eqnarray*}
0
&=& \det
\left[  \matrix {
                 1      & x_{12} & x_{13} & x_{14}  \cr
                 x_{21} & 1      & x_{23} & x_{24}  \cr
                 x_{31} & x_{32} & 1      & x_{34}  \cr
                 x_{41} & x_{42} & x_{43} & 1       \cr
} \right].
\end{eqnarray*}
If some $x_{ij} \not= \pm 1$ then these lines $l_1, \ldots, l_4$ 
are unique up to orientation-preserving isometry and simultaneous 
reversal of each line's orientation.
\end{cor}

It was known to Fenchel (see \cite[\S V.3]{Fenchel}) that if four lines 
$l_1, \ldots, l_4 \in \OrLines$ are given then the complex numbers 
$\lacute l_i, l_j \racute$ satisfy the above hextet equation. 
 
We finish this section with a brief discussion of degenerate arrangements of 
lines in $\hyp$.   

\begin{defn}
\label{D:HexDegen}
An arrangement of lines $l_1, \ldots , l_n \in \OrLines$ is {\em 
non-degenerate} if $l_1, \ldots , l_n$ spans $\Lines$
and the arrangement is {\em degenerate} otherwise.  
\end{defn}

The following lemma gives us three more characterisations of degeneracy.  

\begin{lemma}
\label{L:Degen}
Let $l_1, \ldots, l_n \in \OrLines$ be given, let $x_{ij} = \lacute l_i,l_j 
\racute$ for each  $i,j = 1, \ldots, n$ and let $X$ be the $n \times n$ 
matrix whose $(i,j)^{th}$ entry is $x_{ij}$.  If some $x_{ij} \not= \pm 1$ then the 
following are equivalent:
\begin{itemize}
\item[\rm(1)] $l_1, \ldots, l_n$ is a degenerate arrangement of lines.
\item[\rm(2)] $l_1, \ldots, l_n$ have a common perpendicular line.  
\item[\rm(3)] There exists an orientation-preserving isometry taking the 
arrangement $l_1, \ldots, l_n$ to the arrangement $-l_1, \ldots, -l_n$.
\item[\rm(4)] $\rk(X) = 2$.
\end{itemize}
\end{lemma}

\proof
{\bf (1) $\Rightarrow$ (2)}\qua Suppose that $l_1, \ldots, l_n$ are linearly dependent. 
Some $x_{ij} \not= \pm 1$ so without loss of generality we assume  $x_{12} 
\not= \pm 1$ and hence that $l_1$ and $l_2$ have no end-points in common.
Therefore $l_1$ and $l_2$  are linearly independent and they also have a 
common 
perpendicular $n \in \OrLines$ (whose orientation is not unique).  Then for 
each $i = 3, \ldots, n$ there exist $\alpha_1, \alpha_2 \in \mathbb{C}$ 
so that $l_i = \alpha_1 l_1 + \alpha_2 l_2$.  But 
$\lacute l_1, n \racute = \lacute l_2, n \racute = 0$ and so 
$\lacute l_i, n \racute = 0$ and hence $l_i$ is also perpendicular to $n$
(by Lemma \ref{L:InnerProduct}).

{\bf (2) $\Rightarrow$ (3)}\qua If $l_i \in \OrLines$ is perpendicular to $n \in 
\OrLines$ then the half-turn about $n$ takes $l_i$ to $-l_i$.

{\bf (3) $\Rightarrow$ (1)}\qua
Suppose there is some $g \in \isom$ so that $g \cdot l_i = -l_i$ for 
each $i = 1,\ldots, n$.  If $l_1, \ldots, l_n$ span $\Lines$
then $g$ is an isometry which takes every oriented line to the same 
line but with the opposite orientation, which is absurd. 

{\bf (1) $\Rightarrow$ (4)}\qua
If we assume that any three of the lines $l_1, \ldots, l_n$ are linearly 
dependent then it follows that any three of the rows of $X$ are linearly 
dependent, too, and hence that $\rk(X) <3$.  For example, if 
$\alpha_1 l_1 + \alpha_2 l_2 + \alpha_3 l_3 = 0$ then 
\begin{eqnarray*}
0 &=& \alpha_1 \lacute l_1,l_i \racute + \alpha_2 \lacute l_2,l_i \racute 
+ \alpha_3 \lacute l_3,l_i \racute  \cr
&=&  \alpha_1 x_{1i}+\alpha_2 x_{2i} +\alpha_3 x_{3i}
\end{eqnarray*}
for each  $i = 1, \ldots, n$.  But since some $x_{ij} \not= \pm 1$, 
$\rk(X) \not = 1$ and so $\rk(X) = 2$.

{\bf (4) $\Rightarrow$ (1)}\qua
Conversely, suppose that $\rk(X) = 2$ and assume (in order to derive a 
contradiction) that the $l_1, \ldots, l_n$ span $\Lines$.  Without loss 
of generality assume that $l_1, l_2, l_3$ form a basis for $\Lines$.  
Then since $\rk(X) = 2$, the first three rows of $X$ are linearly dependent 
so there exist $\alpha_1, \alpha_2, \alpha_3 \in \mathbb{C}$ so that 
\begin{eqnarray*}
0 &=& \alpha_1 x_{1i}+\alpha_2 x_{2i} +\alpha_3 x_{3i}  \cr
  &=& \lacute \alpha_1 l_1 + \alpha_2 l_2 + \alpha_3 l_3, l_i \racute
\end{eqnarray*}
for each $i = 1, \ldots, n$.  But since $\lacute \cdot, \cdot \racute$ is 
non-degenerate and (by assumption) $l_1, \ldots, l_n$ span $\Lines$, 
this implies $\alpha_1 l_1 + \alpha_2 l_2 + \alpha_3 l_3 = 0$ which is 
a contradiction.  
\finpf

\section{The ortholength invariant}
\label{S:OrthInv}

In this section we define the ortholength invariant of each incomplete 
hyperbolic structure in the hyperbolic Dehn surgery space $\defm$ of $M$.   
This invariant is given purely in terms of holonomy representations 
so its definition naturally extends to a 
map $\orth:\Char \dashrightarrow \mathbb{C}^n$ from the 
$\isom$-character variety $\Char$ of $M$ to $\mathbb{C}^n$.  
We show that $\orth$ is a rational 
map whose image lies inside a variety $\param \subseteq \mathbb{C}^n$.
For examples of $\orth$, $\param$ and $\Char$, see Section \ref{S:Eg}.

Let $M$ be an oriented, finite volume, $1$-cusped hyperbolic $3$-manifold 
(as in Section \ref{S:Intro}) and let $N$ be an embedded horoball 
neighbourhood of the cusp.  
Then $N$ is diffeomorphic to the product of a $2$-torus $T^2$ with the 
half-open interval $[0,1)$ and furthermore 
the complement in $M$ of the interior of $N$ 
is a compact $3$-manifold with boundary $\partial N \cong T^2$ (e.g.\ see 
Thurston \cite[\S 4.5]{ThurstonBook}).  Let $K$ be a (topological) ideal 
triangulation 
of $M$ (see Benedetti-Petronio \cite[\S E.5-i]{BP}) which meets $N$ `nicely', 
i.e.\ so that inside any tetrahedron $\triangle$ of $K$, $N$ has four 
connected components, each being a punctured neighbourhood of one of the 
vertices of $\triangle$.  Let $*$ be a  base-point for $M$
which lies in $\partial N$.

Now, let $\rho: \pi_1(M,*) \to \isom$ be a homomorphism which satisfies the 
condition that  
$\rho(\pi_1(\partial N,*))$ fixes exactly two points on the sphere 
at infinity $\si$ (e.g.\ $\rho$ could be the holonomy 
representation of an incomplete hyperbolic structure of $\defm$).
A simple investigation of the fixed-points of the Abelian 
group $\rho(\pi_1(\partial N,*))$ shows that this condition is equivalent to 
the requirement that $\rho(\pi_1(\partial N,*))$ is a non-trivial
group of non-parabolic\footnote{An orientation-preserving isometry
$g \in \isom$ is {\em parabolic} if $g \not= 1$ and $\tr^2g = 4$.} 
isometries which is not isomorphic to
$\mathbb{Z}/2 \oplus \mathbb{Z}/2$ (where $\mathbb{Z}/2$ is the group
with two elements).  These conditions in turn are equivalent to the 
algebraic conditions that 
\begin{eqnarray}
\label{E:OrthDomain}
\mbox{neither\quad} \tr^2\rho(m) = \tr^2\rho(l) = 4 \cr
\mbox{nor\quad} \tr^2\rho(m) = \tr^2\rho(l) = \tr^2\rho(ml) =  0, 
\end{eqnarray}
where $m$ and $l$ are any pair of generators of 
$\pi_1(\partial N,*) \cong \mathbb{Z} \oplus \mathbb{Z}$.  

Now, by our assumption that $\rho$ satisfies the conditions 
(\ref{E:OrthDomain}), it follows that the group $\rho(\pi_1(\partial N,*))$ 
fixes a unique line in $\hyp$.  Choose an orientation and call the resulting 
oriented line $\sigma$.
Let $\pi: \tM \to M$ be the universal cover of $M$ and choose
some base-point $\tilde{*} \in \pi^{-1}(*)$.  This choice
allows us to identify the deck-transformations of $\pi:
\tM \to M$ with $\pi_1(M,*)$.  Let $\tN$ be the connected component of 
$\pi^{-1}(N)$ which contains $\tilde{*}$.

Let $e_1, \ldots, e_n$ be the edges of $K$ and for each 
$i = 1, \ldots, n$ choose a lift $\tilde{e}_i$ of edge $e_i$ in 
$\tM$.  To each end of $\tilde{e}_i$ there 
is a corresponding connected component of $\pi^{-1}(N)$.  Denote these 
connected components by $\tN_1$ and $\tN_2$ (note that it is possible 
that $\tN_1 =\tN_2$ if $e_i$ is homotopically trivial).  
Then since $\tN_1$ and $\tN_2$ both cover $N$, 
there exist deck-transformations $\gamma_1$ and $\gamma_2$ for which 
$\gamma_i(\tN) = \tN_i$ (for each $i = 1,2$).  We define the {\em ortholength}
$d_i$ corresponding to edge $e_i$ to be 
$$d_i \defeq \dist(\rho(\gamma_1) \cdot \sigma, \rho(\gamma_2) \cdot \sigma)$$
i.e.\ define $d_i$ to be the complex distance between 
$\rho(\gamma_1) \cdot \sigma$ and $\rho(\gamma_2) \cdot \sigma$.  
While $\gamma_1$ and $\gamma_2$ are not unique, 
$\rho(\gamma_1) \cdot \sigma$ and $\rho(\gamma_2) \cdot \sigma$ 
do not depend on their arbitrariness and so are 
well-defined oriented lines in $\hyp$.  
Also, the definition of $d_i$ doesn't depend on our choice of lift 
$\tilde{e_i}$ of edge $e_i$.  This is because any other lift is of the form 
$\alpha \cdot \tilde{e_i}$ for some deck-transformation $\alpha$.  
In the above prescription this has the effect of replacing 
$\tN_i$ by $\alpha(\tN_i)$, i.e.\ replacing $\gamma_i$ by 
$\alpha \gamma_i$ ($i = 1,2$).  But since 
$\dist(\rho(\alpha) \rho(\gamma_1) \cdot \sigma, \rho(\alpha) \rho(\gamma_2) 
\cdot \sigma) =\dist(\rho(\gamma_1) \cdot \sigma, \rho(\gamma_2) \cdot 
\sigma)$, $d_i$ is not affected by this change.  

Also, conjugating $\rho$ by some $g \in \isom$ has the effect of replacing 
$\sigma$ by $g \cdot \sigma$ and each $\rho(\gamma_1)$ by $g\rho(\gamma_1)
g^{-1}$, which clearly leaves $d_i$ unchanged.  This shows that 
$d_i$ is independent of our choice of $\tilde{*}$ (which we used to 
identify $\pi_1(M,*)$ with the deck-transformations of $\widetilde{M} \to M$) 
and also that $d_i$ only depends on the conjugacy class of $\rho$.  
We define the ortholength invariant $\orth(\rho)$ with respect to 
the ideal triangulation $K$ evaluated at $\rho$ to be 
$$\orth(\rho) \defeq (\cosh d_1, \ldots, \cosh d_n) \in \mathbb{C}^n. $$
We think of $\orth(\cdot)$ as the function which associates $\cosh$ of the 
ortholength $d_i$ to edge $e_i$ of $K$ for each $i = 1, \ldots, n$. 

Now, as noted in Footnote \ref{FootnoteRef} (and see also \cite{BZ}) 
$\isom$ is naturally isomorphic 
to $\mbox{SO}_3\mathbb{C}$ and so the space of representations $\mathcal{R}(M)$ 
of $\pi_1(M,*)$ into $\isom$ is a complex algebraic variety.  
The $\isom$-character variety $\Char$ of $M$ is the `quotient' (in the sense 
of algebraic geometry, see \cite[\S\S 3,4]{BZ}) 
of  $\mathcal{R}(M)$ by the conjugacy action of $\isom$.  This space 
$\Char$ has a natural algebraic structure which makes it an affine 
algebraic variety whose co-ordinate ring is the ring of 
regular functions on $\mathcal{R}(M)$ which are invariant under 
the $\isom$-conjugacy action.  

As noted above, $\orth(\cdot)$ is defined on all representations satisfying 
the conditions (\ref{E:OrthDomain})
and $\orth(\cdot)$ is constant on conjugacy classes.    
Hence $\orth(\cdot)$ descends\footnote{Since $\Char$ is not simply the quotient of the representation 
variety $\rep$ 
of $M$ by the conjugacy action of $\isom$ this is not immediately obvious.
We need two more facts: (1) $\orth(\rho) = (1, \ldots, 1) \in 
\mathbb{C}^n$ whenever $\rho$ is reducible (i.e.\ 
whenever $\rho(\pi_1(M,*))$ fixes a point on the sphere at infinity $\si$)
and (2) if $\rho$ is irreducible and $\rho, \rho^\prime \in \rep$ project 
to the same point under the algebro-geometric quotient map $\rep \to \Char$ 
then $\rho$ and $\rho^\prime$ are conjugate (see \cite[p. 753]{BZ}).
} to a function defined on the complement of a 
proper\footnote{The subvariety is proper since it doesn't contain
(the characters of) the holonomy representations of the (incomplete) 
hyperbolic structures of $\defm$.}
sub-variety of $\Char$.  
We denote this function by  $\orth$ as well, and we write 
$\orth:\Char \dashrightarrow \mathbb{C}^n$
to indicate that $\orth$ is not necessarily defined\footnote{%
Warning:  Since $\Char$ is not necessarily irreducible it is 
possible that $\orth$ is not defined at all on some components of $\Char$.} 
as a set-theoretic function at all points of $\Char$. 
 
Now, given a tetrahedron $\triangle$ of $K$, choose some tetrahedron
$\widetilde{\triangle}$ in $\tM$ which covers it.  Each corner of 
$\widetilde{\triangle}$ meets a unique connected component of $\pi^{-1}(N)$ 
and these connected components determine four 
(not necessarily distinct) oriented lines
$\rho(\gamma_1) \cdot \sigma, \ldots, \rho(\gamma_4) \cdot \sigma$ in $\hyp$, 
as above. The complex distance between any pair of these lines is equal to 
the ortholength 
corresponding to one of the edges of $\triangle$.  Hence the four oriented 
lines $\rho(\gamma_1) \cdot \sigma, \ldots, \rho(\gamma_4) \cdot \sigma$
determine a hextet which realises the complex distances associated to the
edges of $\triangle$.  By Corollary \ref{C:N=4} this implies that 
the hyperbolic cosines of these ortholengths satisfy a certain algebraic 
equation for each tetrahedron of $K$.  Hence the image 
of $\orth$ lies in the following (not necessarily irreducible) complex 
algebraic variety $\param \subseteq \mathbb{C}^n$.  

\begin{defn}[The Ortholength Space $\param$]
\label{D:param}
The {\em ortholength space} $\param\! \subseteq \mathbb{C}^n$ corresponding 
to a (topological) ideal triangulation $K$ of $M$ with $n$ edges  
is the complex affine algebraic variety consisting of those points of 
$\mathbb{C}^n$
which satisfy the {\em hextet equations} of all the tetrahedra of $K$.
Here the hextet equation of a tetrahedron  $\triangle$ of $K$ is 
$$0 = \det
\left[  \matrix {
                 1      & x_{12} & x_{13} & x_{14}  \cr
                 x_{21} & 1      & x_{23} & x_{24}  \cr
                 x_{31} & x_{32} & 1      & x_{34}  \cr
                 x_{41} & x_{42} & x_{43} & 1       \cr
} \right]  $$
where the vertices of $\triangle$ have been numbered 
arbitrarily from $1$ to $4$ and where we denote $\cosh$ of the 
ortholength associated to the edge of $\triangle$ between vertices $i$ and $j$ 
by $x_{ij}= x_{ji}$.  
\end{defn}

Note that the hextet equation of $\triangle$
doesn't depend on the arbitrary numbering of its vertices. 

Although $\param$ lies in an $n$-dimensional space (one dimension for each 
edge of $K$) and is defined by $n$ hextet equations (one for each tetrahedron 
of $K$) for `generic' $K$, $\param$ has an irreducible 
component which is a complex curve.   This follows from Theorem \ref{T:Param}
(below) and from the fact that Dehn surgery space is diffeomorphic to 
$\mathbb{C}$ in a neighbourhood of the complete structure (see Thurston 
\cite[\S 5.5]{ThurstonNotes}).  Conversely, it should be possible to 
prove that the irreducible component of $\Char$ which 
contains the complete hyperbolic structure is one-dimensional over $\mathbb{C}$
via Theorem \ref{T:Param} and a lemma about the dimension of $\param$.

We finish this section with an explicit formula for the ortholength invariant. 
This formula is given in terms of a presentation for 
$\pi_1(M,*)$ based on the ideal triangulation $K$ of $M$.  The generators 
for this presentation 
consist of generators for $\pi_1(\partial N,*)$ plus loops $\alpha_i \not\in
\pi_1(\partial N,*)$ which lie in $\partial N \cup e_i$, where $e_i$ is the 
$i^{th}$ edge of $K$ ($i = 1, \ldots, n$).  

\begin{prop}
\label{P:orth}
Let $\rho:\pi_1(M,*) \to \isom$ be a representation on which 
$\orth$ is defined (i.e.\ for which the conditions (\ref{E:OrthDomain}) hold) 
and for $i = 1, \ldots,n$ let $\tilde{h}, \tilde{g}_i \in 
\mbox{SL}_2\mathbb{C}$ be matrices which cover $\rho(\beta), \rho(\alpha_i) 
\in \isom$,  for $\alpha_i$ as above and some non-trivial 
$\beta \in \pi_1(\partial N,*)$.  Then the $i^{th}$ co-ordinate $\cosh d_i$ of 
$\orth(\rho) \in \param \subseteq \mathbb{C}^n$ is
$$\cosh d_i = 2 \frac{ \mbox{tr}( \tilde{h}\tilde{g}_i) \mbox{tr} 
(\tilde{h}^{-1}\tilde{g}_i) - \mbox{tr}^2 \tilde{g}_i   }
         {   \mbox{tr}^2 \tilde{h} - 4   } -1$$
where $i = 1, \ldots, n$.
\end{prop}

Note that $\mbox{tr}( \tilde{h}\tilde{g}_i) \mbox{tr}
(\tilde{h}^{-1}\tilde{g}_i)$, $\mbox{tr}^2 \tilde{g}_i$ and 
$\mbox{tr}^2 \tilde{h}$ are independent of the choice of 
matrices $\tilde{h}$ and $\tilde{g}_i$ covering $\rho(m)$ and $\rho(\alpha_i)$
and so these three functions define elements of the 
co-ordinate ring of $\Char$ (see \cite{BZ}).  Hence from the above formula it 
is clear that $\orth:\Char \dashrightarrow \param$ is a rational map.  

\proof[Proof of Proposition \ref{P:orth}]
Suppose we have the set-up as in the statement, and let
$\rho(\pi_1(\partial N,*))$ fix a geodesic $\sigma$ of $\mathbb{H}^3$, which we
give an arbitrary orientation.  Then $\cosh d_i$ is 
equal to $\cosh(\dist(\sigma, g_i \cdot \sigma))$, where $g_i = \rho(\alpha_i) 
\in \isom$.  Since $\cosh(\dist(\sigma, g_i \cdot\sigma))$ is invariant under
conjugating $\rho$ by an orientation-preserving isometry, we can assume 
\begin{eqnarray*}
\tilde{h} =   \left[ \matrix{ e^{x/2} & 0 \cr 0 & e^{-x/2} \cr} \right] &
\tilde{g}_i = \left[ \matrix{ a & b \cr c & d \cr} \right],
\end{eqnarray*}
(where $ad -bc = 1$) for the purposes of calculating
$\cosh(\dist(\sigma, g_i \cdot \sigma))$.
We'll then express our answer in terms which are invariant under
conjugacy, and then our formula will hold true for general $\tilde{h}$ and 
$\tilde{g}_i$.

The axis of $\tilde{h}$ has end-points $0$ and $\infty$ on the sphere at 
infinity, so the line matrix of $\sigma$ is 
$\pm${\tiny $\left[ \matrix{ i & 0 \cr 0 & -i \cr} \right]$}$
\in \mbox{SL}_2\mathbb{C}$.    Then by Lemma \ref{L:InnerProduct}
and the fact that the line matrix of $g_i \cdot \sigma$ is 
$\pm \tilde{g}_i ${\tiny $\left[ \matrix{ i & 0 \cr 0 & -i \cr} \right] $}$ 
\tilde{g}_i^ {-1}$, we have
\begin{eqnarray*}
\cosh(\dist(\sigma, g_i \cdot \sigma))
&=& - \mbox{$\frac{1}{2}$ tr}
\left[ \matrix{ i & 0 \cr 0 & -i \cr} \right]
\left[ \matrix{ a & b \cr c & d \cr} \right]
\left[ \matrix{ i & 0 \cr 0 & -i \cr} \right]
\left[ \matrix{ d & -b \cr -c & a \cr} \right]   \\
&=& ad + bc \\
&=& 2ad -1.
\end{eqnarray*}
So our task now is to express $ad$ invariantly.  But
\begin{eqnarray*}
\mbox{tr}(\tilde{h}\tilde{g}_i) = e^{x/2}a + e^{-x/2}d & \mbox{ and } &
\mbox{tr}(\tilde{h}^{-1}\tilde{g}_i) = e^{-x/2}a + e^{x/2}d
\end{eqnarray*}
so
\begin{eqnarray*}
\mbox{tr} (\tilde{h}\tilde{g}_i) \mbox{tr} (\tilde{h}^{-1}\tilde{g}_i)
&=&   a^2 + d^2 + ad(e^{x/2} + e^{-x/2}) \\
&=&   (a+d)^2 - 2ad + 2ad \cosh x \\
&=&   \mbox{tr}^2 \tilde{g}_i + 2ad ( \cosh x - 1).
\end{eqnarray*}
Combining this with
$$ \cosh x - 1 = 2 \cosh^2 (x/2) - 2 = (\mbox{tr}^2 \tilde{h}  - 4)/2$$
gives us the required formula for $\cosh d_i$.
\finpf

\section{Parameterising hyperbolic structures with ortho\-lengths}
\label{S:Param}

In this section we prove that the ortholength invariant locally parameterises the deformation space $\defm$ and is a complete invariant when restricted 
to the (topological) Dehn fillings of $M$ which admit a hyperbolic structure 
(as long as the hyperbolic structure in question has ortholength invariant 
not lying in a certain (possibly-empty) subvariety).  We also prove  
(under a weak technical assumption) that the ortholength invariant 
is a birational equivalence from certain interesting irreducible 
components of $\Char$ to certain irreducible components of $\param$.  

Fix a (topological) ideal triangulation $K$ of $M$ with $n$ edges and let a 
set of ortholength parameters  
$p = (p_1, \ldots, p_n) \in \param$ be given.  Then a {\em realisation}
of $p$ is a set of $n$ hextets (see Figure \ref{F:HexTet}) each of which 
realises the ortholength parameters associated to the edges of a tetrahedron 
of $K$.  More precisely, a realisation is an association of an oriented 
line to the four corners of each tetrahedron $\triangle$ of $K$ so 
that $\cosh$ of the complex distance between any pair of these four 
lines is equal to the ortholength parameter of the edge of $\triangle$ which 
lies between the corresponding pair of vertices of $\triangle$.   

The broad aim of this section is to construct an inverse to $\orth:\Char 
\dashrightarrow \param$.  To each point $p \in \param$, Corollary 
\ref{C:N=4} guarantees the existence of a realisation of $p$ by hextets.
Our strategy for building a representation 
$\pi_1(M,*) \to \isom$ is essentially to glue copies of 
these hextets together to form a type of skeletal developing map for $M$ whose 
rigidity then gives us a holonomy representation for free (see the 
proof of Lemma \ref{L:ConRep}, below).  However, this 
approach is complicated slightly by the fact that the ortholength parameters 
$p = (p_1, \ldots, p_n)$ don't quite determine the 
hextets up to isometry (see Corollary \ref{C:N=4}), even if we assume no 
$p_i = \pm 1$.  This forces us to make Definition \ref{D:coherence}, below.    

First note that each face $F$ of $K$ is contained in exactly two tetrahedra 
of $K$ 
and that these two inclusions composed with the realisation of $p$ by 
hextets gives two right-angled hexagons corresponding to $F$.  Here we think 
of a right-angled hexagon corresponding to $F$ as an association of 
an oriented line to each of the three corners of $F$, and by an isometry  
between such hexagons we mean an orientation-preserving isometry which 
respects this association.  

\begin{defn}[Coherence]
\label{D:coherence}
Let $p = (p_1, \ldots, p_n) \in \param$ be such that no $p_i = \pm 1$. Then a realisation of $p$ by hextets is {\em coherent} if the pair of 
right-angled hexagons corresponding to each face of $K$ are isometric.   
If a coherent realisation of $p$ exists then we say that $p$ is {\em coherent}.
\end{defn}

Note that any pair of  right-angled hexagons corresponding to a face of $K$ 
have the same ortholengths and so by Theorem \ref{T:SolnCh1} must be either
isometric or isometric after reversing the orientations of the
lines in one of the hexagons.  Also note (from the definition of ortholength 
invariant, see Section \ref{S:OrthInv}) that any point in 
the set-theoretical image of $\orth$ is coherent.  

Now, since the rigidity of hextets described in Corollary 
\ref{C:N=4} may fail if all 
of the ortholength parameters are $\pm 1$, we will usually restrict our 
attention to $p \in \param \setminus \mathcal{T}$, where  
$$\mathcal{T} = \{(p_1, \ldots, p_n) \in \param \mid  p_i = \pm 1
\mbox{ for some } i = 1, \ldots, n \}.$$

Also, let $F_1, \ldots, F_{2n}$ denote the faces of $K$ and for each $m = 1, 
\ldots, 2n$ define 
$$\mathcal{S}_m \defeq \{ p \in \param \mid 
 \det \left[ \matrix {1 & p_{i} & p_{j} \cr
                           p_{i} & 1 & p_{k} \cr 
                           p_{j} & p_{k} & 1} \right] = 0 \} $$
where $p_i,p_j,p_k$ are the co-ordinates of $p = 
(p_1, \ldots, p_n)$ corresponding to the edges of $F_m$.  By Lemma 
\ref{L:Degen}, $\mathcal{S}_m$ is the set of ortholengths $p\in \param$ 
so that any realisation of $p$ by hextets has a degenerate hexagon 
corresponding to face $F_m$.  Define 
\begin{eqnarray}
\label{E:calS}
\mathcal{S} \defeq \bigcup_{i \not= j} \mathcal{S}_i \cap \mathcal{S}_j
\end{eqnarray}
i.e.\ a point of $\param$ lies in $\mathcal{S}$ only if it has multiple 
degeneracies.  Note that if $\param$ is one 
complex-dimensional (which Theorem \ref{T:Param} indicates whenever $M$ 
has a single cusp) then dimensional considerations suggest $\mathcal{S}$ 
will be empty.  

The following lemma shows how the existence and rigidity of hextets (Corollary 
\ref{C:N=4}) allows us to construct holonomy representations from 
coherent ortholength parameters.  

\begin{lemma}
\label{L:ConRep}
To each {\em coherent} realisation of $p \in \param \setminus \mathcal{T}$ by 
hextets there exists a corresponding 
representation $\rho: \pi_1(M,*) \to \isom$ so that $\orth(\rho) = p$.  
Furthermore, up to conjugacy there are at most a finite number of 
representations $\rho$ for which $\orth(\rho) = p$
and there is only one if $p \not\in \mathcal{S}$.  
\end{lemma}

Hence this lemma implies that, in $\param \setminus \mathcal{T}$, 
the image of $\orth :\Char \dashrightarrow \param$ is exactly the set of 
coherent ortholength parameters.

Before proving Lemma \ref{L:ConRep} we pause briefly to describe part 
of the relationship between Thurston and 
{\em SnapPea's} parameterisation of  $\Char$ via 
ideal hyperbolic tetrahedra (see \cite{ThurstonNotes}, \cite{SnapPea}, 
\cite{NZ}) and the  
parameterisation of  $\Char$ described in Lemma \ref{L:ConRep} via 
ortholengths $p \in \param$.  To each hextet there is an  associated 
ideal hyperbolic tetrahedron whose vertices are the  
second end-points\footnote{The {\em second end-point} of an oriented line 
$(p,q) \in \endpts$ is $q$.} of the four oriented lines comprising the hextet.
So given a coherent realisation of $p \in \param$, we can associate a 
hyperbolic ideal tetrahedron to each (topological) ideal 
tetrahedron of $K$.  Then {\em SnapPea's} edge-consistency 
conditions  are automatically satisfied for these geometric tetrahedra.  
To see this, consider all the tetrahedra of $K$ surrounding an edge $e_i$ of 
$K$.  In a hextet which realises one of these tetrahedra, two of the four 
lines comprising the hextet correspond to edge $e_i$ 
(and hence $\cosh$ of the complex 
distance between these two lines is $p_i$).  By acting by (orientation-preserving) 
isometries if need 
be we can assume that these two lines are the same for each hextet.
But then the pair of right-angled hexagons 
corresponding to a face $F$ of $K$ containing edge $e_i$ 
must also coincide (since they are 
isometric by the assumption of coherence and they share two distinct lines).  
Associating hyperbolic ideal tetrahedra (as above) to each of these hextets  
therefore gives geometrical tetrahedra whose 
faces coincide and which fit together smoothly about their shared edge. 

The shape of a hyperbolic ideal tetrahedron is determined by a complex 
shape parameter (see \cite{ThurstonNotes}).  
The shape parameter $z\in \mathbb{C}$ of the hyperbolic tetrahedron 
associated (by the above procedure) 
to a given hextet satisfies an equation in $\cosh$ of the complex 
distances between the lines of the hextet.  This equation is quite large 
but it is clearly quadratic in 
$z$, reflecting the fact that a set of ortholengths only determine a 
corresponding hextet up to orientation.  

\proof[Proof of Lemma \ref{L:ConRep}]
Let $p \in \param\setminus \mathcal{T}$ be coherent and
let $H_1, \ldots, H_n$ be a coherent realisation of $p$ by hextets, where
$H_j$ corresponds to the $j^{th}$ tetrahedron $\triangle_j$ of $K$.  
We can pull the ideal triangulation $K$ back to an ideal triangulation 
$\widetilde{K}$ of the universal cover $\widetilde{M}$ of $M$ via the 
covering map $\widetilde{M} \to M$.  Then every tetrahedron, face and edge of 
$\widetilde{K}$ covers one of the tetrahedra, faces or edges (respectively) 
of $K$.  Let $\widetilde{K}^{(3)}$ denote the set of tetrahedra 
of $\widetilde{K}$.  We will consider a function $H:\widetilde{K}^{(3)}
\to (\OrLines)^4$ which associates a hextet to each tetrahedron of 
$\widetilde{K}$.  This is as an association of an oriented line 
to each of the four corners of each tetrahedron of $\widetilde{K}$.  
We require that $H$ satisfies two conditions: (1) if $\widetilde{\triangle}
\in \widetilde{K}^{(3)}$ covers the $j^{th}$ tetrahedron $\triangle_j$ of 
$K$ then the hextet
$H(\widetilde{\triangle})$ is isometric to $H_j$ and (2) if two tetrahedra
$\widetilde{\triangle}, \widetilde{\triangle}^\prime \in \widetilde{K}^{(3)}$
share a face $T$ of $\widetilde{K}$ then the right-angled hexagons corresponding 
to $T$ in
$H(\widetilde{\triangle})$ and $H(\widetilde{\triangle}^\prime)$ are 
identical.  

We will show that to each coherent realisation of $p$ by hextets 
there is a corresponding 
function $H:\widetilde{K}^{(3)} \to (\OrLines)^4$ which satisfies the above 
two conditions, and that this function is unique up to (orientation-preserving) 
isometry.  Before proving 
this we first show how the existence and rigidity of such a 
`skeletal developing map' $H$ proves the lemma.  

Given such an $H$, by its uniqueness up to isometry we know that 
for each deck-transformation $\alpha$ of the covering $\widetilde{M} \to 
M$ there is a unique isometry $\rho(\alpha) \in \isom$ which makes the 
following diagram
$$
\begin{array}{ccc}
\widetilde{K}^{(3)} & \stackrel{\alpha}{\to} & \widetilde{K}^{(3)}  \\
\downarrow H & & \downarrow H        \\
(\OrLines)^4 & \stackrel{{\rho(\alpha)}}{\rightarrow} & (\OrLines)^4
\end{array}
$$
commute.  Choosing a base-point $\tilde{*} \in \widetilde{M}$ covering 
$* \in M$ allows us to identify the deck transformations of the covering  
$\widetilde{M} \to M$ with $\pi_1(M,*)$, and hence $\rho$ becomes a 
function $\rho:\pi_1(M,*) \to \isom$.  It follows easily from the definition 
that $\rho:\pi_1(M,*) \to \isom$ is actually a homomorphism, and 
it is well-defined up to conjugacy since using a different 
$H$ or a different base-point $\tilde{*}$ simply has the effect of 
conjugating $\rho$ by an isometry.  

By construction, translates under $\rho(\pi_1(M,*))$ of the line which is 
fixed by $\rho(\pi_1(N,*))$ are 
simply the lines occurring in the image of $H$, so $\orth(\rho) = p$. 
Up to isometry there are at most $2^n$ different realisations of $p$ and 
so there are at most this many representations $\rho$ with 
$\orth(\rho) = p$ (up to conjugacy).   
Also, if  $p \not\in \mathcal{S}$ then up to isometry 
there are at most
two coherent realisations of $p$, and these are the same except that they 
have opposite orientations.  (This is because the orientation of one hextet 
determines the orientations of all other hextets in a coherent realisation of 
$p$ whenever all or all but one of the hexagonal faces are non-degenerate, 
see point (3) of Lemma \ref{L:Degen}.)  These two realisations therefore 
give rise to functions $H:\widetilde{K}^{(3)} \to (\OrLines)^4$ which are  
identical (up to isometry) except that all of their lines are of opposite 
orientation.  Clearly the holonomy representations corresponding to  
such $H$ are identical.

The rest of the proof is devoted to showing that, corresponding to a 
coherent realisation $H_1, \ldots, H_n$ of $p$, 
there is a function $H:\widetilde{K}^{(3)} \to (\OrLines)^4$ which 
satisfies the two conditions listed in the first paragraph of this proof and 
furthermore that such an $H$ is unique up to (orientation-preserving) isometry.

Dual to $\widetilde{K}$ in $\widetilde{M}$ there is a 
$2$-dimensional CW-complex $C$ whose underlying space 
$|C|$ is a deformation retract of $\widetilde{M}$ and so in particular $|C|$ 
is simply connected.  Let
$\gamma$ denote a finite sequence $\gamma_1, \ldots, \gamma_m$ of oriented 
$1$-cells in $C$ so that $\gamma_i(1) = 
\gamma_{i+1}(0)$ for each $i = 1, \ldots, m-1$, where  $\gamma_i(0)$ 
denotes the start of the $1$-cell $\gamma_i$ and $\gamma_i(1)$ denotes the end 
of it.  Let $x_0 = \gamma_1(0)$ and let $x_i = \gamma_i(1)$ for each 
$i = 1, \ldots, m$.  We think of $\gamma$ as a path in the $1$-skeleton of 
$C$ from $x_0$ to $x_m$. 
 
For each $i = 1, \ldots, m$, the $0$-cell $x_i$ is dual to a tetrahedron 
$\widetilde{\triangle}_i$ (say) 
of $\widetilde{K}$, and each tetrahedron $\widetilde{\triangle}_i$ 
covers a tetrahedron $\triangle_{j(i)}$ of $K$ (for some $j(i) \in \{1, 
\ldots, n \}$).  
Given a choice of hextet $H^\prime_0$ isometric to $H_{j(0)}$, 
we define a sequence $H^\prime_0, \ldots, H^\prime_m$ of hextets 
with each $H^\prime_i$ isometric to $H_{j(i)}$ as follows.  Assume that 
$H^\prime_{i-1}$ is defined.  Then the $1$-cell
$\gamma_i$ is dual to some $2$-simplex $\widetilde{T}$ of $\widetilde{K}$ 
which is contained in both $\widetilde{\triangle}_{i-1}$ and 
$\widetilde{\triangle}_i$.  Let $T$ be the $2$-simplex of $K$ covered by 
$\widetilde{T}$.  By the coherence of $p$ we know 
that the right-angled hexagons in $H_{j({i-1})}$ and $H_{j(i)}$ corresponding 
to $T$ are isometric. Since  $H^\prime_{i-1}$ is isometric to 
$H_{j(i-1)}$ we can define $H^\prime_{i}$ as the hextet which is isometric to 
$H_{j(i)}$ and for which the right-angled hexagons in  $H^\prime_{i-1}$ and 
$H^\prime_i$ corresponding to $\widetilde{T}$ are 
identical.  This uniquely determines $H^\prime_i$ because 
there is no non-trivial (orientation-preserving) isometry fixing 
the right-angled hexagon
corresponding to $\widetilde{T}$, since no $p_i = \pm 1$ and so the hexagon 
has at least three end-points on $\si$.  Hence to any such path $\gamma$
and any choice $H^\prime_0$ we have a corresponding sequence 
$H^\prime_0, \ldots, H^\prime_m$ of hextets.

Now, choose  some $\widetilde{\triangle} \in \widetilde{K}^{(3)}$
(which covers the $j^{th}$ tetrahedron $\triangle_j$ of $K$, say) 
and choose a hextet  which is isometric to $H_j$ and denote it by 
$H(\widetilde{\triangle})$. 
Let $*$ be the $0$-cell of $C$ dual to  $\widetilde{\triangle}$.  Given 
any other $\widetilde{\triangle}^\prime \in \widetilde{K}^{(3)}$, 
dual to some $0$-cell $x$ of $C$, choose a path $\gamma$ as above 
with $x_0 = *$ and $x_m = x$.  Then form the corresponding sequence 
$H^\prime_0, \ldots, H^\prime_m$ of hextets as defined above with 
$H^\prime_0 = H(\widetilde{\triangle})$ 
and define $H(\widetilde{\triangle}^\prime)$ to be $H^\prime_m$.  

To see that this definition  is 
independent of the path from $*$ to $x$, note that any two such paths 
are homotopic in $|C|$.  It therefore suffices to show that the holonomy 
around the boundary of any of the $2$-cells of $C$ is trivial, i.e.\ that 
if $\gamma$ is a path tracing once around a $2$-cell of $C$ then in the 
corresponding sequence $H^\prime_0, \ldots, H^\prime_m$ of hextets defined 
above, $H^\prime_0 = H^\prime_m$.  But for such a $\gamma$, all of the 
tetrahedra $\widetilde{\triangle}_0, \ldots, \widetilde{\triangle}_m$ share 
an edge in $\widetilde{K}$ and so by construction 
all of the hextets $H^\prime_0, \ldots, H^\prime_m$
have two lines in common.  So the hextets $H^\prime_0$ and $H^\prime_m$ 
share two lines and are isometric (since $\widetilde{\triangle}_0 = 
\widetilde{\triangle}_m$).  Then since there 
are no non-trivial orientation-preserving isometries fixing two 
distinct\footnote{Since no 
$p_i = \pm 1$, the two lines have more than three end-points on $\si$.
} lines in $\hyp$, $H^\prime_0 = H^\prime_m$ as required.  

It is clear that 
$H:\widetilde{K}^{(3)} \to (\OrLines)^4$ as constructed satisfies 
the two defining conditions given in the first paragraph of this proof,  
and conversely any such $H$ can be constructed in this way.  Hence the 
existence and rigidity of $H$ is proved. 
\finpf

Non-coherent ortholengths correspond to 
representations into $\isom$ of the fundamental groups of finite-sheeted 
covers or branched covers of $M$, branched over the edges of $K$.  Hence
non-coherent ortholengths seem to be more related to the ideal 
triangulation $K$ than to intrinsic properties of the $3$-manifold $M$.

We next show that the property of coherence is 
locally constant on the complement of a certain subvariety of $\param$.  
This proof is a deformation argument in which we make use of the fact that 
in some sense there exists a
continuous local parameterisation of hextets by their ortholengths.  We pause 
now briefly to give a precise description of this local parameterisation.  

Define the set of {\em reduced hextets} $\redh$ to be 
$$\redh \defeq \{ (l_1, \ldots, l_4) \in (\OrLines)^4 \mid 
l_1\! =  \left[ \matrix {i & 0\cr 0 & -i\cr} \right]\! , 
l_2\! =  \left[ \matrix {a & i-a\cr i+a & -a\cr} \right]\!, 
 a \not= \pm i \in \mathbb{C} \}, $$
i.e.\ the set of hextets whose lines $l_1, \ldots, l_4$ 
are in  standard position so that 
$l_1$ has ordered end-points $(0,\infty)$ and $l_2$ has ordered end-points 
$(b,1)$ for some $b \in \si$ not equal to $0$, $1$  or $\infty$.
Note that any hextet $l_1, \ldots, l_4$ with $\lacute l_1,l_2 \racute 
\not= \pm 1$ is isometric to exactly one of the hextets in $\redh$.  Now, let 
$$\Pd = \{ (x_{12},x_{13},x_{14},x_{23},x_{24},x_{34}) \in 
\mathbb{C}^6 \mid x_{12} \not= \pm1 \mbox{ and }\det X = 0 \} $$
where 
\begin{eqnarray}
\label{E:4X}
X \defeq \left[ \matrix
{1 & x_{12} & x_{13}&x_{14} \cr x_{12} & 1 & x_{23} & x_{24} \cr
x_{13} & x_{23} & 1 & x_{34} \cr x_{14} & x_{24} & x_{34} & 1 \cr } \right]
\end{eqnarray}
and let $\Od:\redh \to \Pd$ be the quadratic map given in the above 
co-ordinates by $x_{ij} = \lacute l_i, l_j \racute$.  By Lemma 
\ref{L:InnerProduct}, $\Od$ takes a hextet and gives us $\cosh$ of the 
ortholengths between the hextet's four lines.  

Now, by Corollary \ref{C:N=4} and Lemma  \ref{L:Degen}, $\Od$ is onto and 
it is two-to-one everywhere except on the set of degenerate hextets where 
$\Od$ is instead one-to-one.  By Lemma  \ref{L:Degen} again, the 
degenerate hextets are the pre-image of those elements of $\Pd$ for which 
$\rk (X) = 2$.   The restriction of $\Od$ 
to the non-degenerate hextets is a double cover onto the elements of 
$\Pd$ for which $\rk (X) = 3$.   
The deck-transformation for this cover is the involution $\sigma: \redh \to 
\redh$ given by 
$\sigma :(l_1,l_2,l_3, l_4) \mapsto (l_1, l_2, f_a(l_3), f_a(l_4))$ where 
$$ f_a( \left[ \matrix {u & v\cr w & -u\cr} \right]) = 
\left[ \matrix {u & \frac{i-a}{i+a}w\cr 
\frac{i+a}{i-a}v & -u\cr} \right]$$
and $a \not= \pm i$ is determined by  
$l_2 =   \left[ \matrix {a & i-a\cr i+a & -a\cr} \right]$.

Similarly, we define the set of reduced right-angled hexagons to be 
$$\rhex \defeq \{ (l_1, l_2, l_3) \in (\OrLines)^3 \mid 
l_1 =  \left[ \matrix {i & 0\cr 0 & -i\cr} \right] , 
l_2 =  \left[ \matrix {a & i-a\cr i+a & -a\cr} \right], 
 a \not= \pm i \in \mathbb{C} \}, $$
and define 
$$\mathcal{P}(F) = \{ (x_{12},x_{13},x_{23}) \in \mathbb{C}^3 
\mid x_{12} \not= \pm1 \} $$
and let $\OrthF:\rhex \to \mathcal{P}(F)$ be the map given in the above 
co-ordinates by $x_{ij} = \lacute l_i, l_j \racute$. 
As above, $\OrthF$ restricts to a double covering from the set of 
non-degenerate right-angled hexagons to the subset of $\mathcal{P}(F)$ 
for which 
$$ \det \left[ \matrix {1 & x_{12} & x_{13} \cr
                           x_{12} & 1 & x_{23} \cr 
                           x_{13} & x_{23} & 1} \right] \not= 0.  $$
  
\begin{lemma}
\label{L:Coher1}
The property of coherence is locally constant on 
$\param \setminus (\mathcal{S} \cup \mathcal{T}).$
\end{lemma}

\proof
Let $a$ be a point of $\param \setminus (\mathcal{S} \cup 
\mathcal{T})$ and for each $i = 1, \ldots, n$, choose 
a numbering from $1$ to $4$ for the vertices of the $i^{th}$ 
tetrahedron $\triangle_i$ of $K$.  Then we have a projection   
$\pi_i : \param \to \mathcal{P}(\triangle)$ essentially given by 
ignoring the edge parameters which do not appear on the edges of $\triangle_i$.
We know that $\pi_i(a)$ is a smooth point of the hypersurface 
$\mathcal{P}(\triangle) \subseteq \mathbb{C}^6$, since $\Od$ is a 
local diffeomorphism onto a neighbourhood of $\pi_i(a)$ (since $\pi_i(a)$ 
corresponds to a non-degenerate hextet) and $\redh$ is easily seen to be 
smooth everywhere.
Let $V_i\subseteq \mathcal{P}(\triangle)$ be a neighbourhood of 
$\pi_i(a)$ which is diffeomorphic to $\mathbb{C}^{5}$ and  
which is disjoint from the hypersurface given by the equation $\rk(X) = 2$, 
where $X$ is as defined in (\ref{E:4X}).  (This is possible since 
$a \not\in \mathcal{S}$ and hence $\pi_i(a)$ corresponds to a 
non-degenerate hextet.)  Let $U$ be a 
connected neighbourhood of $a$ in $\param \setminus (\mathcal{S} 
\cup \mathcal{T})$ which is contained in $\pi_i^{-1}(V_i)$
for each $i = 1, \ldots, n$.  

Now, suppose that there exists some $p\in U$ which is coherent, say with a 
coherent realisation $H^p_1, \ldots, H^p_n \in \redh$, and let 
$q$ be any other point of $U$.  For each $i = 1, \ldots, n$, the image of $U$ 
under $\pi_i : \param \to \mathcal{P}(\triangle)$ is contained in the 
path-connected, locally path-connected and simply connected set $V_i$.  Also, 
the restriction of $\Od:\redh \to \mathcal{P}(\triangle)$ to 
$\Od^{-1}(V_i)$ is a covering projection.  Hence by the lifting 
property of coverings (see \cite{Bredon}) there is a unique lift 
$\tilde{\pi}_i: U \to \redh$ of  the restriction $\pi_i|U$ 
of $\pi_i$ to $U$ for which $\tilde{\pi}_i(p) = H^p_i$.  We define $H^q_i$ 
to be the hextet $\tilde{\pi}_i(q)$.  

Now, let $F$ be any face of $K$ and suppose that $F$ is contained in 
tetrahedra $\triangle_i$ and $\triangle_j$ of $K$.  
If we choose a numbering from $1$ to $3$ for 
the vertices of $F$ then we can define several projection maps.  Firstly 
a projection $\pi_F:\param \to \mathcal{P}(F)$ essentially given by forgetting 
all ortholength parameters except those associated to the edges of $F$.  
Also, a map $\tau_i: \redh \to \rhex$ which first  
forgets the line of $H^q_i\in \redh$ corresponding to 
the vertex of $\triangle_i$ not contained in $F$ and which then 
acts by an (orientation-preserving) isometry to bring the resulting 
right-angled hexagon into standard 
position.  Similarly we have a map  $\tau_j: \redh \to \rhex$.

Now, if the right-angled hexagons 
corresponding to $F$ in $H^q_i$ and $H^q_j$ are degenerate then 
they are isometric by Lemma \ref{L:Degen}.  Hence we can restrict 
our attention to the set $U^\prime$ consisting of $q \in U$ 
so that $H^q_i$ and $H^q_j$ are non-degenerate.  
Then $\tau_i \circ \tilde{\pi}_i: U^\prime \to \rhex$ and 
$\tau_j \circ \tilde{\pi}_j: U^\prime \to \rhex$ are both lifts of 
the restriction $\pi_F|U^\prime$ of $\pi_F$ to $ U^\prime $ and 
they agree at $p$.  Hence by the uniqueness of lifts to covering 
spaces (see \cite{Bredon}), $\tau_i \circ \tilde{\pi}_i  
= \tau_j \circ \tilde{\pi}_j$ and so the right-angled hexagons 
corresponding to $F$ in $H^q_i$ and $H^q_j$ are again isometric.  
\finpf

\begin{lemma}
\label{L:Coher2}
The limit of a sequence of coherent points in $\param$ is coherent.  
\end{lemma}

\proof
This lemma follows from the fact that the inverse images of 
compact sets in $\Pd$ under $\Od : \redh \to \Pd$ are compact. 
\finpf

A rational map between two irreducible complex varieties which is 
generically\footnote{%
A property is {\em generically} satisfied on a variety if it is true on the 
complement of a proper subvariety.} defined, 
injective and onto has a rational inverse and is called a birational 
equivalence (see \cite[p.77]{Harris}). Hence  
Lemmas \ref{L:ConRep}, \ref{L:Coher1} and \ref{L:Coher2} combine to give  
the following theorem.  

\begin{theorem}
\label{T:Param}
Let $C$ be an irreducible component of $\Char$ on which $\orth$ is 
generically defined and suppose that the set-theoretic image of $C$ under 
the rational map $\orth:\Char \dashrightarrow \param$ contains some point 
$p = (p_1, \ldots, p_n)\in \param$ for which  
$p \not\in \mathcal{S}$ (see equation (\ref{E:calS})) and no $p_i = \pm 1$.  
Then $\orth$ restricts to a birational equivalence 
between $C$ and an irreducible component of $\param$.  
\end{theorem}  
 
Let $\hol : \defm \to \Char$ be the map which takes a hyperbolic structure 
and returns its holonomy representation modulo conjugacy.  Then as in 
Section \ref{S:OrthInv}, the ortholength invariant of an incomplete 
hyperbolic structure $M_0 \in \defm$ is $\orth \circ \hol (M_0)$.  
Recall also that we say that an incomplete hyperbolic structure on $M$ {\em corresponds 
to a (topological) Dehn filling} if the metric completion of $M$ 
is a closed manifold (homeomorphic to a topological Dehn filling of $M$) 
with a smooth hyperbolic structure.  

\begin{theorem}
\label{T:CP}
Let $M_0 \in \defm$ be an incomplete hyperbolic structure whose ortholength 
invariant $p = (p_1, \ldots, p_n) \in \param$ is not 
contained in $\mathcal{S}$ (see equation (\ref{E:calS})) and which has no 
$p_i = \pm 1$.  Then:
\begin{itemize}
\item If $M_0$ corresponds to a (topological) Dehn filling of $M$ then 
$M_0$ is uniquely determined by $p$.  
\item The ortholength invariant $\orth : \Char \dashrightarrow \param$ 
restricts to a diffeomorphism from a neighbourhood of $\hol(M_0)$ in 
$\Char$ to a neighbourhood of $p$ in $\param$.
\end{itemize}
\end{theorem}

The holonomy map $\hol$ is a local homeomorphism (e.g.\ see 
\cite[\S 5.2]{ThurstonNotes} or \cite{Goldman}) so we can endow $\defm$  
with a smooth structure coming from $\Char$.  Then the second point of this 
theorem says that (generically) the ortholength invariant smoothly 
locally parameterises incomplete hyperbolic structures on $M$.  

\proof[Proof of Theorem \ref{T:CP}]
A hyperbolic structure which corresponds to a topological Dehn filling of 
$M$ is determined by its holonomy representation.  This is a `folk-lore'
result which is true essentially because, in this case, the 
image of the holonomy representation is a discrete Kleinian group, and its
action on $\hyp$ has a quotient isometric to the metric completion 
of $M$.  But by Lemma \ref{L:ConRep}, the ortholength invariant $p$ 
determines the holonomy representation whenever $p \not\in 
\mathcal{S}$ and no $p_i = \pm 1$.  
This proves the first point of the corollary.  

Now, $\orth$ is defined on every representation which arises as the 
holonomy of one of the hyperbolic structures of $\defm$ (see Section 
\ref{S:OrthInv}).  Then since the conditions (\ref{E:OrthDomain}) are open,
$\orth$ is defined in a neighbourhood of $\hol(M_0)$.  
From the formula of Proposition \ref{P:orth} it is clear that 
$\orth$ is smooth in this neighbourhood.

On the other hand, $p$ is coherent by the construction used to 
define the ortholength invariant in Section \ref{S:OrthInv}.  
By Lemma \ref{L:Coher1} this implies that all ortholength parameters 
in a neighbourhood of $p$ are also coherent.  Hence by Lemma 
\ref{L:ConRep} there is a local inverse to $\orth$ defined in a 
neighbourhood of $p$.  This inverse is smooth by the construction 
given in Lemma \ref{L:ConRep}, since the representation 
given there is determined by a finite portion of the skeletal 
developing map $H$ and the hextets used to build this map 
all vary smoothly with $p$.  This last result follows from 
the discussion preceding Lemma \ref{L:Coher1} and the fact that 
if $p \not\in \mathcal{S}$ then the hextets in a realisation of $p$ 
are all non-degenerate.  
\finpf

Theorems \ref{T:Param} and \ref{T:CP} are vacuous unless there exists 
some incomplete hyperbolic structure with ortholength 
invariant $p$ such that $p \not\in \mathcal{S}$ and no $p_i = \pm 1$.  
This will probably be the case whenever the edges of 
$K$ are homotopically non-trivial, however until such a general result 
can be established we have the following lemma.  

\begin{lemma}
\label{L:EP}
If the ideal triangulation $K$ is Epstein-Penner's \cite{EP} ideal 
cell-decomposition of $M$ then there exists 
an incomplete hyperbolic structure with ortholength 
invariant $p$ such that $p \not\in \mathcal{S}$ and no $p_i = \pm 1$.

\end{lemma}

The full proof of this lemma relies on the notion of a `tube domain' 
(see Section \ref{S:Conc} and \cite[Lemma 4.7]{Dowty}) so we  
content ourselves here with a sketch.  

\proof[Sketch of a proof of Lemma \ref{L:EP}]  
Corresponding to the complete hyper-\break bolic structure on $M$ there is a
2-complex in $M$ which is the image of the boundary of the Ford 
domain of $M$ under the face-pairing identifications.  This {\em Ford 2-complex} 
is dual to $K$ (see \cite{EP}).  Let $\mbar$ be a  
(topological) Dehn filling of $M$ with very large Dehn surgery co-ordinates 
and let $\Sigma$ be the added `core' geodesic of $\mbar$, hence $\mbar$ is 
hyperbolic and $M = \mbar \setminus \Sigma$.  
The cut-locus of $\Sigma$ in $\mbar$ is a 2-complex which (for very large Dehn 
fillings) is a small perturbation of the Ford 2-complex, and 
so in particular this cut-locus is also topologically dual to $K$ (see \cite{Dowty}).  
Now, the pre-image 
of $\Sigma$ under a (fixed) covering isometry $\hyp \to \mbar$ is a collection 
of disjoint lines in $\hyp$.  If $F_j$ is a 2-simplex of $K$ and 
$\widetilde{F}_j$ is a lift to $\hyp$ then each corner of $\widetilde{F}_j$ 
determines one of these lines.  The complex distances between these three 
distinguished lines are given by the three ortholength parameters $p_i$ attached to the 
corresponding edges of $F_j$, where $p = (p_1, \ldots, p_n)$ is the ortholength 
invariant of the hyperbolic structure on $M$ induced from $M \hookrightarrow \mbar$.  
Then the set of points of $\hyp$ which are equidistant from all three lines 
simultaneously is non-empty, since this set contains a lift of a certain part 
of the cut-locus of $\Sigma$, namely the 1-cell dual to $F_j$.

But if $p$ lies in $\mathcal{S}_j$ then this is a contradiction.  
In this case, the three lines corresponding to the 
corners of $\widetilde{F}_j$ have a common perpendicular line $n$ (by Lemma 
\ref{L:Degen}).  Since $\mbar$ is a large Dehn filling, all of its ortholengths 
will be very large\footnote{For large Dehn fillings the tube radius is large, 
see \cite{ThurstonNotes}.} and so the distances between 
these three lines are also large (and no $p_i = \pm 1$).  
This implies that the locus of points equidistant from any pair of 
lines is close to 
the plane which is perpendicular to $n$ and lies mid-way between the two lines.  
Hence there can be no point of $\hyp$ which is equidistant from all 
three lines simultaneously.  
\finpf

\section{Examples}
\label{S:Eg}

In this section we describe the ortholength invariant $\orth: \Char 
\dashrightarrow \param$ for the figure-8 knot complement $M$.

\begin{figure}[ht!]
\begin{center}
\includegraphics[height=3cm]{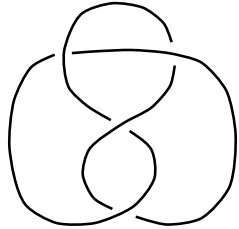}
\end{center}
\caption{The figure-8 knot}
\label{F:Fig8}
\end{figure}

\begin{figure}[ht!]
\vspace{-0.5cm}
\begin{center}
\subfigure[Tetrahedron $0$]{
\includegraphics[height=4cm]{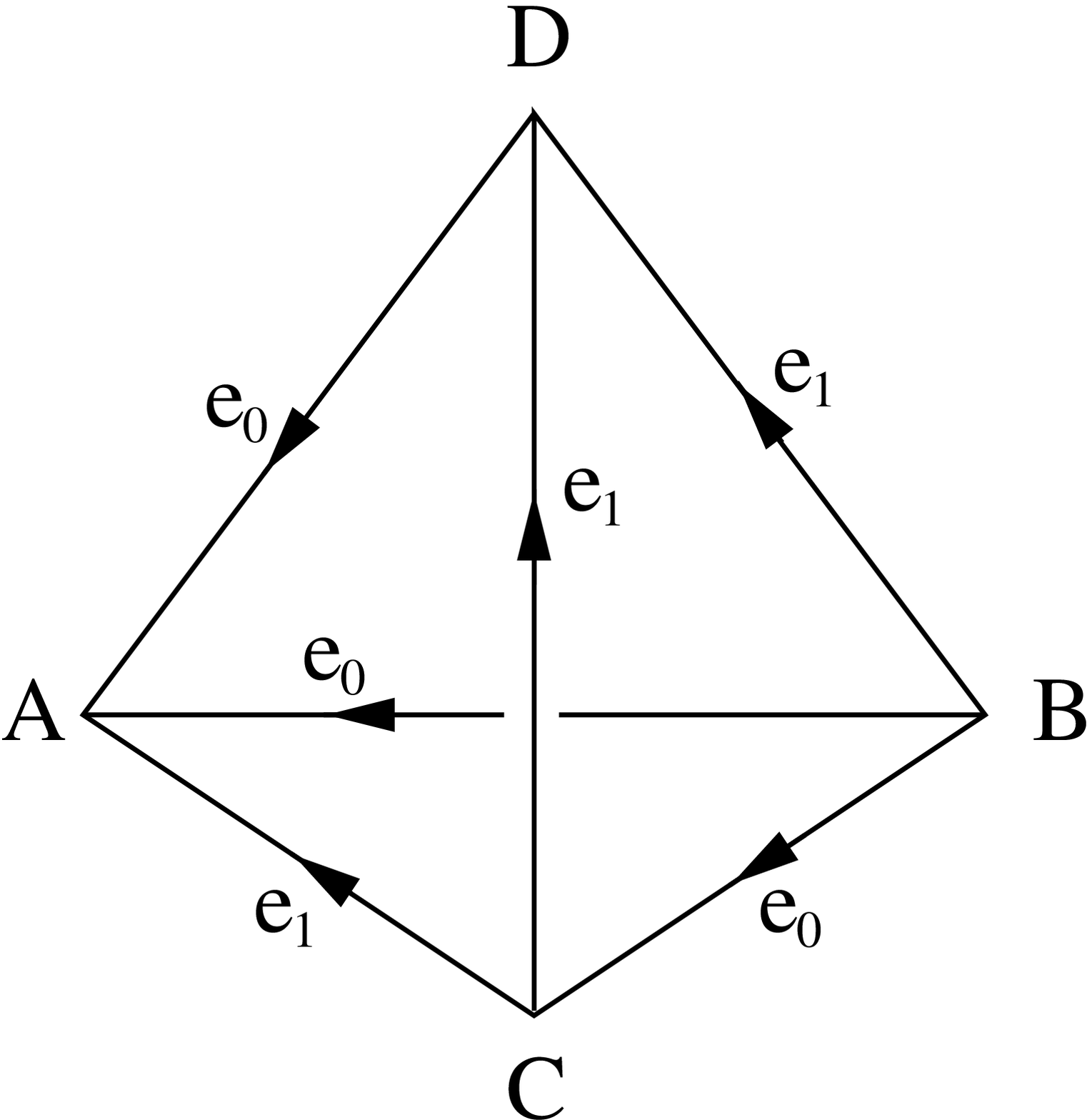}\hspace{1cm}
}
\subfigure[Tetrahedron $1$]{
\includegraphics[height=4cm]{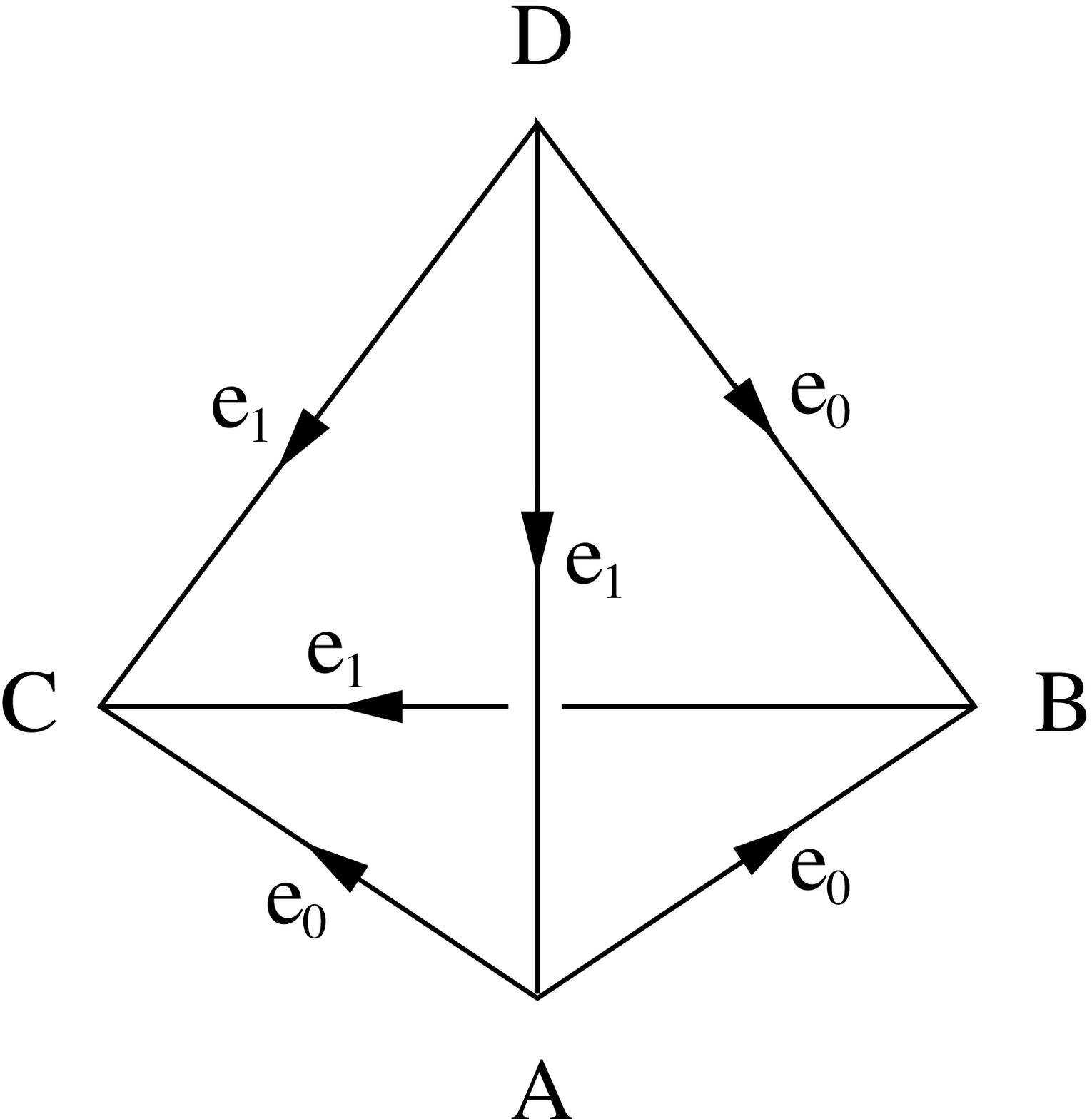}
} 
\end{center}\vspace{-0.5cm}
\caption{The canonical ideal triangulation $K$ for the figure-8 knot complement}
\label{F:tetn_FP}
\end{figure}

Let $M$ be the complement in $\mathbb{S}^3$ of the figure-8 knot
(see Figure \ref{F:Fig8}).  
A result of Epstein and Penner (see \cite{EP}) says that every 
1-cusped hyperbolic 3-manifold has a canonical ideal cell-decomposition.
For the figure-8 knot complement $M$, this cell-decomposition is 
the ideal triangulation $K$ shown\footnote{The 
faces of the tetrahedra shown in Figure \ref{F:tetn_FP} are identified in 
pairs giving face-classes A,B,C and D. The class of each face is written at 
the vertex opposite the face.} in Figure \ref{F:tetn_FP} (see 
\cite{ThurstonBook}).  This ideal triangulation for the figure-8 knot 
complement was first described by Thurston \cite{ThurstonNotes}.  

The ideal triangulation $K$ has two edge classes, $e_0$ and $e_1$.  
Let $p_0$ be the ortholength parameter corresponding to edge 
class $e_0$ and let $p_1$ correspond to $e_1$.   
Then the hextet equation (see Definition \ref{D:param}) of tetrahedron $0$ is 
\begin{eqnarray}
0 &=& \det  \left[ 
\matrix {1   & p_0 & p_1 & p_1 \cr 
         p_0 & 1   & p_1 & p_0 \cr
         p_1 & p_1 & 1   & p_0 \cr
         p_1 & p_0 & p_0 & 1    } 
\right] \nonumber\\  
&=&  1-3p_0^2 -3 p_1^2 + 4 p_0 p_1^2 + 4 p_0^2 p_1  
+ p_0^4 - 2p_0^3 p_1 - p_0^2 p_1^2 - 2p_0 p_1^3 + p_1^4 \nonumber \\
&=& \mbox{$\frac{1}{4}$} (p_0^2 + p_1^2 + p_0 p_1 - p_0 -p_1 - 1) 
(2p_1 - 3 p_0+1+\sqrt{5} (p_0-1))                            \nonumber \\
&& ~~\times (2p_1 - 3p_0 +1-\sqrt{5} (p_0-1))               \label{E:ParamEqn}.
\end{eqnarray}
The hextet equation of tetrahedron $1$ is identical, 
so  $\param \subseteq \mathbb{C}^2$ is the plane algebraic curve given 
by equation (\ref{E:ParamEqn}).  From (\ref{E:ParamEqn}) 
we see that $\param$ is the union of a conic 
(whose projective completion is topologically a smooth sphere) 
and two lines.  

\begin{figure}[ht!]
\vspace{-0.3cm}
\begin{center}
\subfigure[The generators for\newline $\pi_1(\partial N,*)$]{
\includegraphics[height=3.5cm]{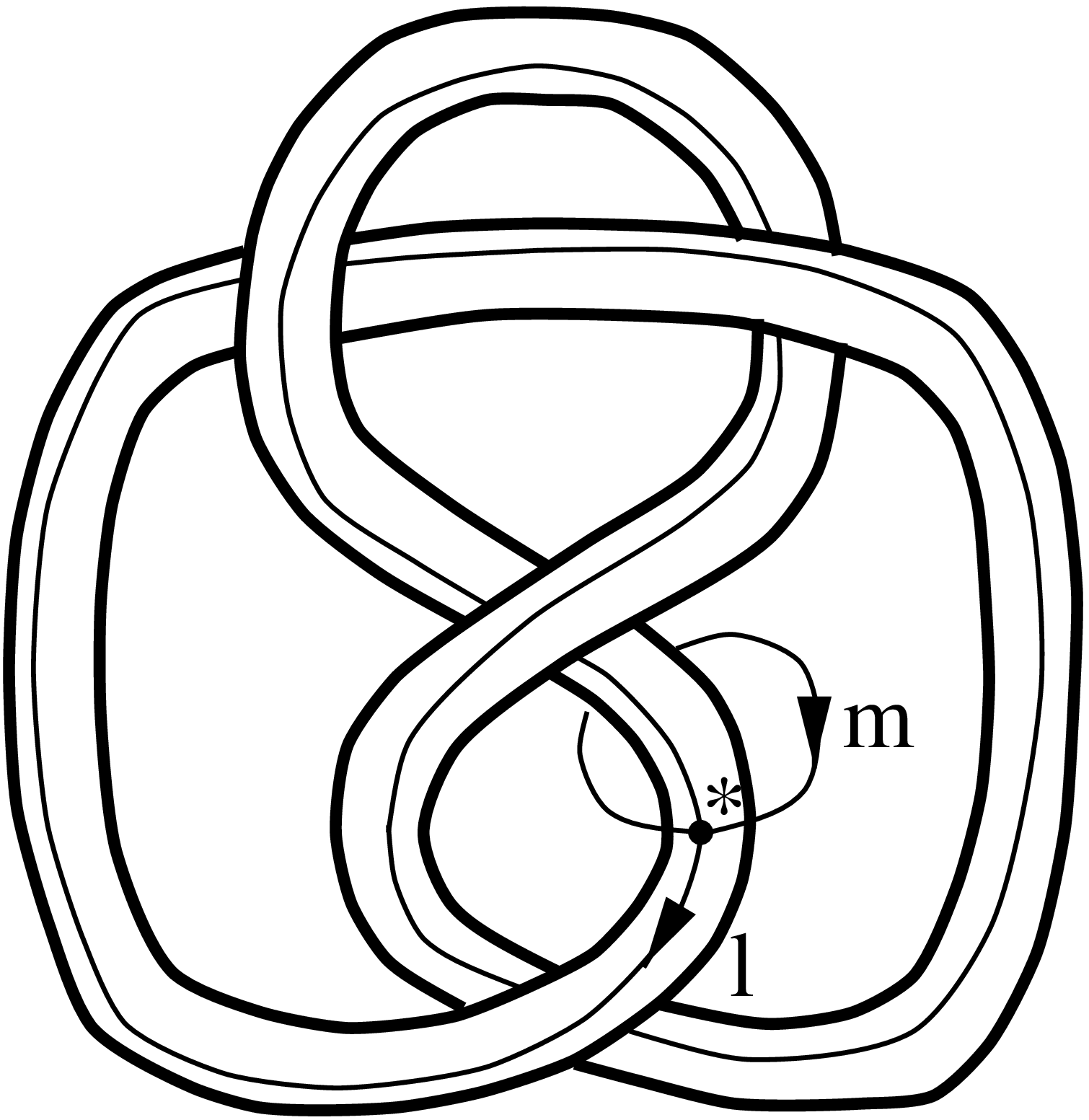}\hspace{0.5cm}
}
\subfigure[The generator $\alpha_0$]{
\includegraphics[height=3.5cm]{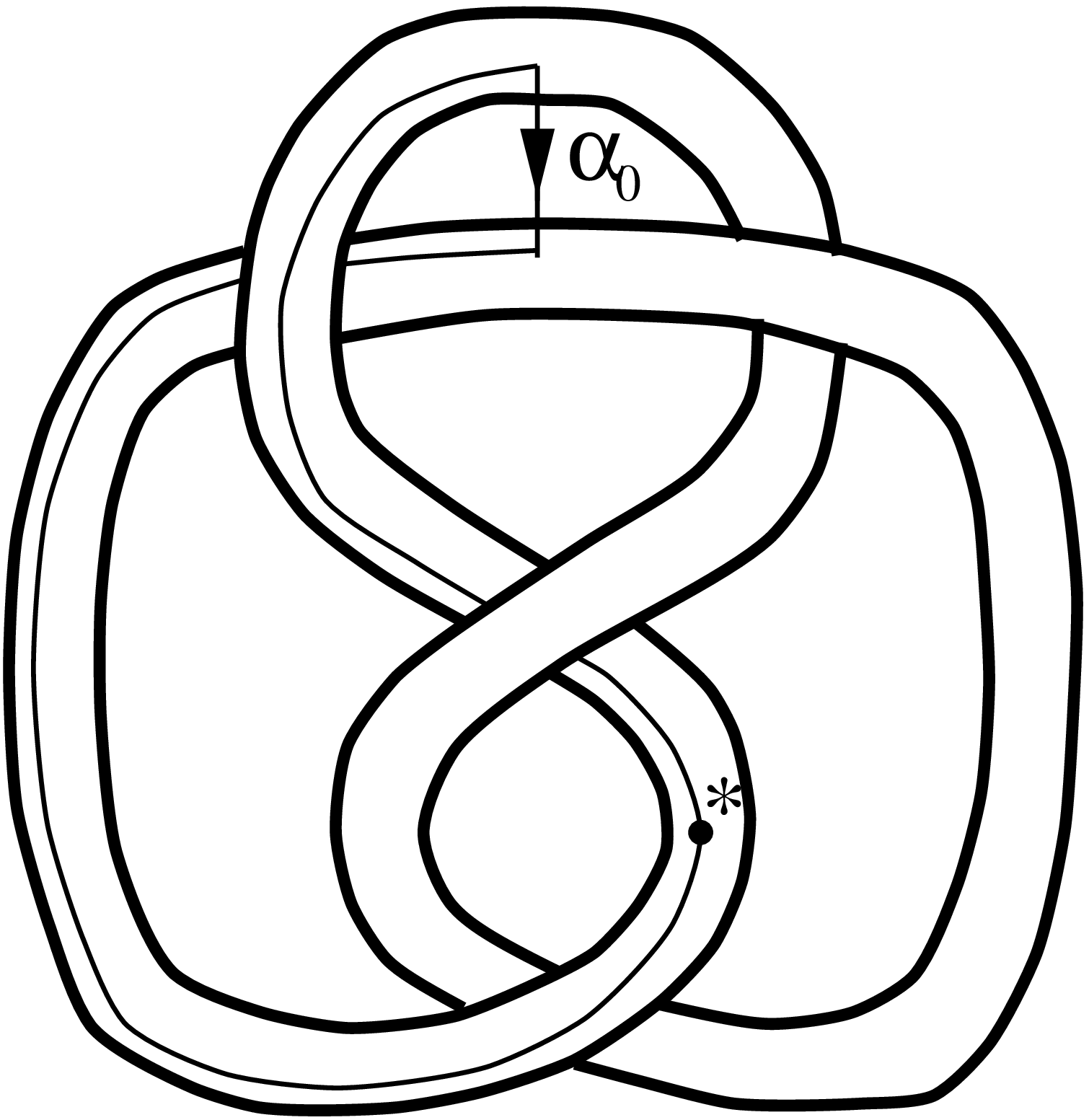}\hspace{0.5cm}
}
\subfigure[The generator $\alpha_1$]{
\includegraphics[height=3.5cm]{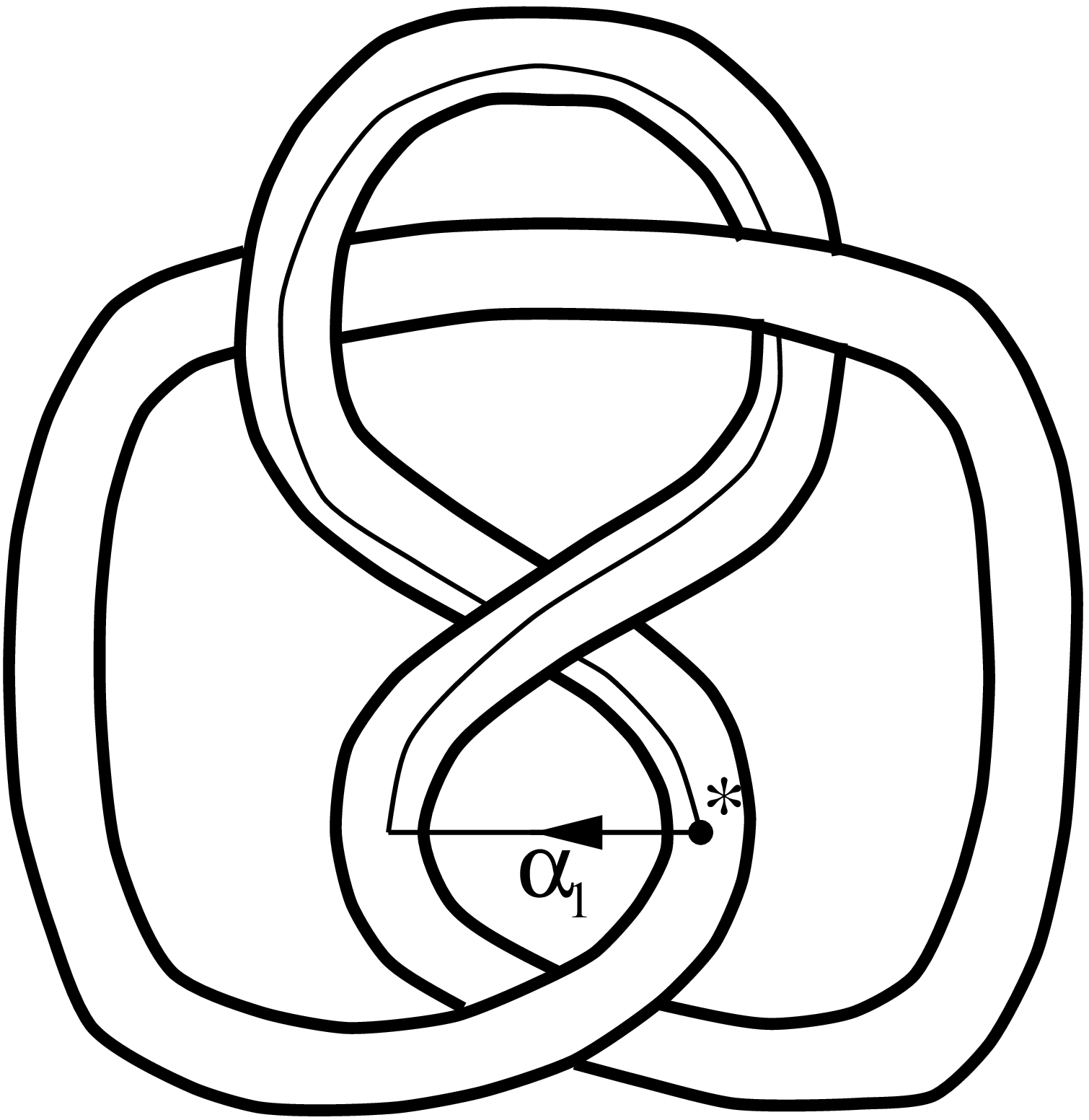}
} 
\end{center}
\vspace{-0.5cm}
\caption{The loops corresponding to the generators of the presentation 
for $\pi_1(M,*)$ based on $K$}
\label{F:Fig8_alphan}
\end{figure}

Following the notation of Section \ref{S:OrthInv} we take $N \subseteq M$ to 
be a closed tubular neighbourhood of the figure-8 knot minus the knot itself.  
Then the generators of the presentation of $\pi_1(M,*)$ based on $K$ 
(see the discussion preceding Proposition \ref{P:orth}) are as shown in Figure 
\ref{F:Fig8_alphan}.  

Now, let $t_1 = l \alpha_1^{-1} m$ and $t_2 = m$.  These are the generators of 
a presentation 
\begin{eqnarray}
\pi_1(M,*) = \lacute t_1, t_2 | [t_2^{-1}, t_1] t_2^{-1} = 
t_1^{-1} [t_2^{-1}, t_1] \racute                            \label{E:WirtPres}
\end{eqnarray}
which is derived from a Wirtinger presentation for $\pi_1(M,*)$.
We define two complex algebraic functions $X$ and $Y$ on the space of 
representations of $\pi_1(\!M,\!*)$ into $\mbox{SL}_2\mathbb{C}$ by 
$$X(\tilde{\rho}) = \mbox{tr}(\tilde{\rho} (t_1)) = 
\mbox{tr}(\tilde{\rho}(t_2)) \mbox{ and } Y(\tilde{\rho}) = 
\mbox{tr}(\tilde{\rho} (t_1 t_2))$$
for any representation $\tilde{\rho}: \pi_1(M,*) \to \mbox{SL}_2 \mathbb{C}$.
These functions are clearly conjugacy invariant, so they descend to functions 
(also denoted $X$ and $Y$) on the $\mbox{SL}_2\mathbb{C}$-character 
variety\footnote{Similarly to the $\isom$-character variety, the 
$\mbox{SL}_2\mathbb{C}$-character variety of $M$ is defined to be the 
algebro-geometric quotient of the space of representations $\pi_1(M,*) \to 
\mbox{SL}_2 \mathbb{C}$ by the conjugacy action of $\mbox{SL}_2\mathbb{C}$,
as defined in Culler-Shalen \cite{CS}.}.  
In fact,  $X$ and $Y$ are the co-ordinate functions of an embedding 
of the $\mbox{SL}_2 \mathbb{C}$-character variety in $\mathbb{C}^2$ as the 
complex curve given by 
\begin{eqnarray}
\label{E:SL2CEqn}
0 = (X^2 - Y -2) (X^2 Y - 2 X^2 -Y^2 +Y +1 )
\end{eqnarray}
(see Gonz\'alez-Acu\~na and Montesinos-Amilibia's paper \cite{Montesinos}).  

Now, from the presentation (\ref{E:WirtPres}) above it is clear that 
each representation $\rho: \pi_1(M,*) \to \isom$ lifts to exactly two 
representation $\pm \tilde{\rho}: \pi_1(M,*) \to \mbox{SL}_2 \mathbb{C}$.   
Hence there are well-defined algebraic functions 
$U$ and $V$ defined on the\break $\isom$-character variety $\Char$ by
$$U(\rho) = X^2(\tilde{\rho}) \mbox{  and  } V(\rho) = X^2(\tilde{\rho}) 
-Y(\tilde{\rho}),$$
i.e.\ $U$ and $V$ are elements of the co-ordinate ring of $\Char$.
By (\ref{E:SL2CEqn}), these functions define an embedding of $\Char$ into 
$\mathbb{C}^2$ as the complex curve given by the equation 
\begin{eqnarray}
0 = (V-2)(V^2 - UV + V + U -1).    \label{E:CharUVEqn}
\end{eqnarray}
Hence $\Char$ is composed of two irreducible components, one being the 
line $V = 2$ and the other being the rational curve
\begin{eqnarray}
U = \frac{V^2 + V -1}{V - 1}.    \label{E:CharParam}
\end{eqnarray}
The line $V=2$ corresponds to reducible representations (i.e.\ those which 
fix a point on $\si$) and the rational curve  (\ref{E:CharParam}) 
corresponds to conjugacy classes of irreducible representations 
(see \cite{Montesinos}).

We can now calculate $\orth:\Char \dashrightarrow \param$ using the formula 
given in Proposition  \ref{P:orth}.  The conditions (\ref{E:OrthDomain}) are 
satisfied at all points of $\Char$ except the points with $(U,V)$ equal to 
$(4,2)$ or $(4,1 - \zeta)$, where $\zeta \in \mathbb{C}$ is a non-trivial 
cube root of unity.  (These points correspond to the trivial representation 
and to the two orientation-preserving conjugacy classes of 
discrete and faithful representations.)  
Since the line $V=2$ corresponds to reducible representations 
(see \cite{Montesinos}) $\orth$ takes the whole line to the point $(1,1)$.
On the curve (\ref{E:CharParam}),  Proposition  \ref{P:orth} implies that
\begin{eqnarray}
\orth(U,V) = \frac{1}{V^2 - 3V +3}(-V^2 + 3V -1, V^2-V -1).
\label{E:lambda}
\end{eqnarray}  From this expression it is easy to check that the image of $\orth$ 
lies inside the conic component of $\param$.  From (\ref{E:lambda}) 
and Figure \ref{F:tetn_FP}, the ortholength invariant 
$\orth(U,V)$ lies in $\mathcal{S}$ whenever 
$$ 0 = (V-2)^2(V^2 +V-1) \mbox{ and } 0=(V -2)^2(V-1)^2(V^2+V-1),$$
where $(U,V)\in \Char$ belongs to the curve  (\ref{E:CharParam}),   
i.e.\ $\orth(U,V) \in \mathcal{S}$ when $V$ equals $2$ or 
$(-1 \pm \sqrt 5)/2$.  Also, a simple calculation shows that 
$\orth(U,V) \in \param \subseteq \mathbb{C}^2$ has one or both 
co-ordinates equal to $\pm 1$ if and only if $V$ is $1$ or $2$.   
Hence by Theorem \ref{T:Param}, the image of $\orth$ is dense in the conic 
component of $\param$ and $\orth$ restricts to a birational equivalence
between this component and the curve (\ref{E:CharParam}).  
A simple calculation shows that the inverse of $\orth$ on this curve 
is given by $(p_0,p_1) \mapsto (U,V)$, where 
$$V = \frac{2p_0 + p_1 +1}{p_0+1}$$ 
and $U$ is given by the equation (\ref{E:CharParam}).

\section{Future applications}
\label{S:Conc}

The work presented in this paper arose out of an attempt to 
construct families of incomplete hyperbolic structures on a cusped 
$3$-manifold $M$ and to estimate the boundary of the deformation 
space $\defm$ in Dehn surgery co-ordinates.

The connection between the ortholength invariant and these two 
problems is the concept of a 
{\em tube domain} (see \cite{Dowty}) which comes from collaborative work 
with Craig D. Hodgson.   A tube domain can be defined for any   
(incomplete) hyperbolic structure in $\defm$, but the description is simplest 
for those hyperbolic structures which correspond to topological 
Dehn fillings of $M$.  If $\mbar$ is a hyperbolic Dehn filling of $M$ then 
there is a covering isometry $\pi:\hyp \to \mbar$ and a distinguished simple 
closed geodesic $\Sigma$ in $\mbar$ so that $\mbar \setminus \Sigma$ is 
diffeomorphic to $M$. We define the {\em tube domain} $\dom$ of the 
Dehn filling to be (the closure of) 
the set $\widetilde{\dom}$ of points of $\hyp$ which are closer to a line 
$\widetilde{\Sigma} \subseteq \pi^{-1}(\Sigma)$ than to any other line
of $\pi^{-1}(\Sigma)$, modulo the deck-transformations which preserve
$\widetilde{\Sigma}$. There is a surjective 
map $\dom \to \mbar$ which makes the diagram 
$$
\begin{array}{ccc}
\widetilde{\dom} & \hookrightarrow & \hyp  \\
\downarrow & & ~\downarrow \pi       \\
\dom & \rightarrow & \mbar
\end{array}
$$
commute.  The tube domain $\dom$ is diffeomorphic to a solid torus and its 
boundary is broken into (non-planar) faces which the map $\dom \to \mbar$ 
isometrically identifies in pairs.  

On the other hand, the ortholength invariant of the hyperbolic structure 
$\mbar \setminus \Sigma$ on $M$ determines the set of lines $\pi^{-1}(\Sigma)$ 
up to isometry (as detailed in the proof of Lemma \ref{L:ConRep}).  
These in turn determine the corresponding tube domain $\dom$ 
and its face-pairing isometries, and so give us back the original hyperbolic 
structure $\mbar \setminus \Sigma$ on $M$.  This gives us a method of 
constructing hyperbolic structures via tube domains which can be generalised to 
calculate hyperbolic structures with more general Dehn surgery-type 
singularities, including hyperbolic cone-manifolds whose singular sets are 
simple closed geodesics.  This approach has been automated by
Goodman-Hodgson in a computer program \cite{tube} called {\em Tube}.
{\em Tube} naturally lends itself to calculating the tube radius and it 
has been used to calculate hyperbolic structures
which {\em SnapPea} \cite{SnapPea} fails to compute.
For instance, {\em Tube} has calculated a hyperbolic 
structure for m004(4,0) (the hyperbolic cone-manifold with cone-angles 
$\pi/2$ along the figure-8 knot) and for m007(3,1) (the manifold 
of third lowest volume on {\em SnapPea's} closed census) while  {\em SnapPea} 
fails to compute a hyperbolic structure for m004(4,0) 
and there is no known positively oriented tetrahedral decomposition of
m007(3,1).  Also, the cone-manifolds 
over the figure-8 knot complement with cone angle 
strictly less than $2\pi/3$ all have hyperbolic structures (see \cite{HLM} 
and also \cite[\S 29]{HodgsonThesis}) 
but for cone angles bigger than $2\pi/4.767 \ldots$ {\em SnapPea's} 
tetrahedral construction fails.  In fact, Choi \cite{Young} 
has shown that there cannot exist a positively-oriented hyperbolic 
ideal triangulation for the cone-manifold over the figure-8 knot complement 
when cone-angles are $2 \pi / 4.767 \ldots$ (i.e.\ when one of the 
tetrahedra of the canonical ideal triangulation flattens out).  

The second problem which motivates the study of the ortholength invariant 
is the estimation of the boundary of $\defm$.  
The ortholength invariant seems especially suited to this task
because of its close relationship with the tube radius and 
because of some suggestive connections between the tube 
radius and the degeneration of hyperbolic structures on $M$ (e.g.\ see  
Kojima \cite{KojPi} and Hodgson-Kerckhoff \cite{Kerckhoff}, \cite{HK}).
For a given (incomplete) hyperbolic structure $M_0 \in \defm$, there is a 
(topological) ideal triangulation $K_0$ so that the tube radius of $M_0$ is 
given in terms of $\orthz(M_0)$.  In particular, the tube 
radius of $M_0$ will be non-zero whenever all co-ordinates of $\orthz(M_0)$ 
lie in $\mathbb C \setminus [-1, 1]$.  

Each edge of this special ideal triangulation $K_0$ corresponds to a pair of 
faces of the tube domain of $M_0$.  The combinatorics 
of the tube domains and hence of the corresponding ideal triangulations $K_0$ are 
locally constant in $\defm$ and also at the complete structure (see 
\cite[\S 4.5]{Dowty}).  Since $\defm$ is non-compact near its boundary, it is 
conceivable that infinitely many inequivalent triangulations $K_0$ may arise 
as $M_0$ varies over $\defm$.  However, we conjecture that this in not the case, 
i.e.\ for any given $M$ only a finite number $K_1, \ldots, K_m$ of (topological) 
ideal triangulations are needed to compute the tube radii of the incomplete 
hyperbolic structures in $\defm$.  
For our irreducible manifold $M$, the conjecture of \cite{Kerckhoff} 
reduces to the assertion that a family of hyperbolic structures degenerates if 
and only if their tube radii approach zero.  This motivates the following conjecture.  

\begin{conj}
\label{C:Conjecture}
Let $M$ be the underlying smooth manifold of an orientable, $1$-cusped hyperbolic 
3-manifold of finite volume.  Then there are a finite number of (topological) ideal 
triangulations $K_1, \ldots, K_m$ of $M$ so 
that a sequence of (incomplete) hyperbolic structures in $\defm$ has a limit in 
$\defm$ if and only if 
the corresponding $m$ sequences of ortholength invariants 
with respect to the ideal triangulations $K_1, \ldots, K_m$ 
all have limits with co-ordinates $\cosh d_i$ in $\mathbb{C} \setminus [-1, 1]$. 
\end{conj}

For instance, when $M$ is the figure-8 
knot complement, computational evidence suggests that this conjecture is true 
with $m = 1$ and $K_1$ the ideal triangulation of Section \ref{S:Eg}.  
More specifically, in Dehn surgery co-ordinates the zero-volume set 
(which in this case is conjectured to be the boundary of $\defm$, see 
\cite{ThurstonNotes}, \cite{HodgsonThesis}) appears to coincide with the 
frontier of the set of hyperbolic structures with $p_0, p_1 \not\in 
[-1, 1]\subseteq \mathbb{C}$, where $(p_0, p_1)$ is  
the ortholength invariant of Section \ref{S:Eg}.  

Practical conditions to determine when the tube domain construction 
fails may provide rigorous estimates about the boundary of $\defm$, 
especially in special cases such as the figure-8 knot complement.  
Analogies between the construction (in practice) of Dirichlet domains and 
the construction of tube domains suggest possible necessary conditions for 
the failure of the tube domain construction though no useful estimates have 
emerged from this approach so far. 

We can also attempt to use {\em ortho-angles}
(which are dual to the ortholengths) to parameterise $\defm$.
The ortho-angles are more directly related to the Dehn surgery co-ordinates
of $\defm$ than the ortholengths and they are very well-behaved near the
complete structure.  This makes the ortho-angle parameterisation of $\defm$
quite a promising line of research.  

Finally we note that even though we have restricted our attention in this 
paper to single-cusped $3$-manifolds, the main definitions and 
proofs\footnote{In particular, Theorems \ref{T:Param} and \ref{T:CP} 
are true in the case of multiple cusps.} 
trivially extend to $3$-manifolds with multiple cusps.

\newcommand{\noopsort}[1]{} 
\newcommand{\printfirst}[2]{#1}
\newcommand{\singleletter}[1]{#1} 
\newcommand{\switchargs}[2]{#2#1}

\Addresses\recd

\end{document}